\documentclass[11pt]{article}
\usepackage[margin=1in]{geometry}
\usepackage[T1]{fontenc}
\usepackage{lmodern}

\newcommand{\mono}[1]{{\fontfamily{pcr}\selectfont #1}}
\usepackage[authoryear,round]{natbib}
\usepackage{amsmath,amssymb,amsthm}
\usepackage{graphicx}
\usepackage{booktabs}
\usepackage{array}
\usepackage{xcolor}

\usepackage{caption}
\usepackage{hyperref}
\hypersetup{colorlinks=true, linkcolor=blue!50!black, citecolor=blue!50!black, urlcolor=blue!50!black}

\newtheorem{theorem}{Theorem}
\newtheorem{lemma}{Lemma}
\theoremstyle{remark}
\newtheorem{remark}{Remark}


\newcommand{\doi}[1]{\href{https://doi.org/#1}{doi:\nolinkurl{#1}}}
\newcommand{\nuG}{\nu_{\Gamma}}
\newcommand{\nuh}{\hat\nu}
\newcommand{\nuA}{\nu_{\Gamma\mathrm{A}}}
\newcommand{\nuB}{\nu_{\Gamma\mathrm{B}}}
\newcommand{\nuC}{\nu_{\Gamma\mathrm{C}}}

\newcommand{\Ainf}{A_\infty}

\title{Analytic and optimized approximations\\to the median of the gamma distribution}
\author{Richard F.~Lyon\\[2pt]
\small Google Research, Mountain View, California, USA\\
\small \texttt{dicklyon@acm.org}}
\date{\today}

\begin{document}
\maketitle

\begin{abstract}
The median $\nuG(\alpha)$ of the gamma distribution has no closed form, but its behavior at both ends is known:
it is asymptotic to $2^{-1/\alpha}e^{-\gamma}$ (with $\gamma$ being the Euler--Mascheroni constant) as the shape parameter
$\alpha \to 0$, and approaches $\alpha - \frac{1}{3}$ as $\alpha \to \infty$.  We write approximations as interpolations
$(1 - g)\,\nu_0 + g\,\nu_\infty$ between two elementary-function bounds or approximations, each exact in one limit,
with an interpolator $g(\alpha)$ going from 0 to 1, and compare three such forms:  a first one between previously
proven bounds, a second ending in an approximation from the beta median that is exact at $\alpha = 1$ but does not
quite approach $\alpha - \frac{1}{3}$ as $\alpha \to \infty$, and a third one, between bounds proved here, that corrects
the large-$\alpha$ end of the second.  For each form, the ideal interpolator's slope at 0, approach to 1, and value at
$\alpha = 1$ are derived in closed form.
Rational interpolators with every coefficient set in closed form from these properties, with no numerical
fitting, reach a worst relative error of $0.0011$ over all $\alpha > 0$, with absolute error below $0.0002$.
Fitting some of these or other coefficients numerically to minimize the worst absolute-or-relative error
then trades accuracy against complexity:  the Pareto frontier reaches $0.000087$ with two fitted constants
and $0.000041$ with five, and from two on uses degree-4 rationals of the second form with no $\alpha^2$ term in the
numerator.  The proofs that the end functions of the
first and third forms are bounds are machine-checked in Lean~4.

\end{abstract}

\section{Introduction}
\label{sec:intro}

The median $\nuG(\alpha)$ of the standard gamma distribution with shape parameter $\alpha > 0$ is the solution
$x$ of $P(\alpha, x) = \frac{1}{2}$, where $P(\alpha, x) = \gamma(\alpha, x)/\Gamma(\alpha)$ is the regularized
lower incomplete gamma function.  It has no closed form in elementary functions except at $\alpha = 1$, where the
distribution is exponential and $\nuG(1) = \log 2$; at $\alpha = 2$ and $\alpha = \frac{1}{2}$, it can be written with the Lambert $W$
function and the inverse error function (Section~\ref{sec:limits}).  The median has been studied mostly for large $\alpha$:
\citet{chen1986bounds} bounded it between $\alpha - \frac{1}{3}$ and $\alpha$,
\citet{choi1994medians} and \citet{berg2006chen} found the leading terms of its asymptotic series, and \citet{berg2006chen, berg2008convexity}, \citet{adell2008ramanujan},
\citet{you2017approximation}, and \citet{chen2017median} extended these results.  \citet{berg2006chen} also found the
small-$\alpha$ asymptote $2^{-1/\alpha}e^{-\gamma}$, where $\gamma$ is the Euler--Mascheroni constant.

In earlier work \citep{lyon2021closed}, the author wrote the median as $\nuG(\alpha) = 2^{-1/\alpha}(A(\alpha) + \alpha)$,
conjectured that $A(\alpha)$ lies between its limits $e^{-\gamma}$ (as $\alpha \to 0$) and $\log 2 - \frac{1}{3}$ (as
$\alpha \to \infty$), which was later proved \citep{lyon2023tight}, and approximated $A(\alpha)$ by interpolating
between those limits with a one-parameter arctangent or rational function.  The stronger conjecture that $A(\alpha)$
decreases monotonically \citep{lyon2021closed, lyon2023tight} remains open.

A study of the beta median \citep{lyon2026median} led to a different approximation.  A ratio form
$g(\alpha)/(g(\alpha) + g(\beta))$ for the median of the beta distribution with shape parameters $\alpha$ and $\beta$ is
exact on the line $\alpha = \beta$ for any $g$, and it is exact on the lines $\alpha = 1$ and $\beta = 1$ if and only if
$g(\alpha) \propto 1/(2^{1/\alpha} - 1)$ (Lemma~4 there).  Scaling this $g$ to the exact gamma median $\log 2$ at $\alpha = 1$
gives, with the name $\nu_{\Gamma 4}$ from \citet{lyon2026median},
\[
\nu_{\Gamma 4}(\alpha) = \frac{\log 2}{2^{1/\alpha} - 1} \, ,
\]
which \citet{lyon2026median} proved is an upper bound on $\nuG$ for $\alpha < 1$ and a lower bound for $\alpha > 1$ (Theorem~3
there).  It has the right form at both ends, but not the right coefficients:  as $\alpha \to 0$ it is $2^{-1/\alpha}\log 2$,
about $23\%$ high, and as $\alpha \to \infty$ it approaches $\alpha - \frac{1}{2}\log 2$, about $0.0132$ below $\alpha - \frac{1}{3}$.

Here we combine and extend these ideas.  We treat each approximation as an interpolation between two
elementary functions, one correct as $\alpha \to 0$ and one as $\alpha \to \infty$, with an interpolator
$g(\alpha)$ that runs from 0 to 1.  We study three such pairs, which we call forms:  Form A, the form of
\citet{lyon2021closed}, between the proven bounds; Form B, ending in $\nu_{\Gamma 4}$, which gives the most accurate fits; and Form C, new here, which corrects
Form B's denominator so that its large-$\alpha$ end is exact to $O(1)$, and whose two ends we prove are lower and upper
bounds.  When the two functions are a lower and an upper bound, as in Form A or Form C,
the interpolation is a convex combination of them, with weight $g$ between 0 and 1.  Form B is the exception:  for
$\alpha > 1$ both of its functions lie below the median, and its interpolator must exceed 1 there.  Each form is a single
formula, one smooth elementary-function expression for all $\alpha > 0$.

\begin{figure}[tbp]
\centering
\includegraphics[width=0.9\linewidth]{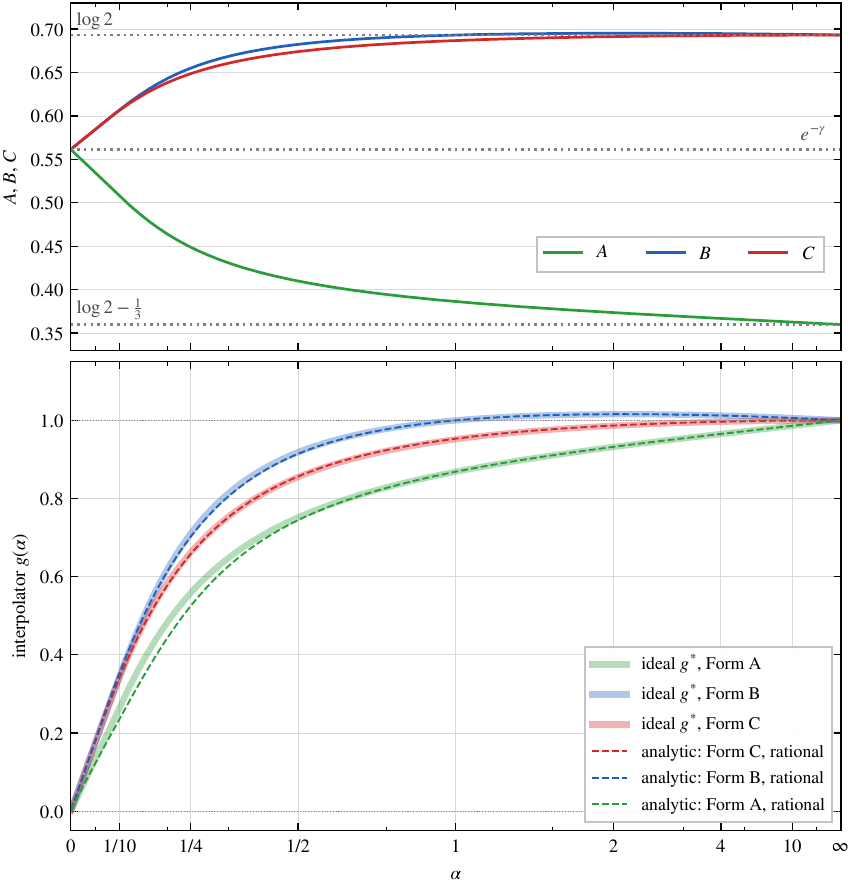}
\caption{Top:  the functions $A$, $B$, and $C$ of \eqref{eq:ABC}, with their limits (dotted).  All three tend to $e^{-\gamma}$
as $\alpha \to 0$; $A$ tends to $\log 2 - \frac{1}{3}$ and $B$ and $C$ to $\log 2$ as $\alpha \to \infty$, $B$ from above.
Bottom:  the ideal interpolators $g^*$ of the three forms (thick, light), which are the same functions rescaled to run from 0 to 1
(Section~\ref{sec:ideal}), with the best analytic rational-function interpolator of each form (Section~\ref{sec:analytic}),
every coefficient set from $g^*$, overlaid (thin, dashed).
Interpolators with fitted coefficients would be indistinguishable from the ideal at this scale.
The horizontal axis is linear in $\arctan\alpha$, here and in subsequent figures that plot functions of $\alpha$.}
\label{fig:ABC}
\end{figure}

For each form, the ideal interpolator $g^*$ that makes it exact is a simple function of the median, and the bounds
above are statements about its range.  Of course, there is no explicit formula for $g^*$, as there is none for the
median itself, but we can compute it numerically and plot it to get some insight into the ideal behavior:  Figure~\ref{fig:ABC}
(bottom) shows the three, and Section~\ref{sec:ideal} discusses their shapes.

We optimize the approximations in two stages.  The first uses no error criterion.  The slope of $g^*$ at $\alpha = 0$, its
approach to 1 at large $\alpha$, and its value at $\alpha = 1$, where the median is $\log 2$, are known in closed form, and a
rational interpolator can take its coefficients from them, as \citet{lyon2021closed} did for Form A\@.  With every coefficient
set this way, the best approximation, of Form C, has a worst relative error of $0.0011$ over all $\alpha > 0$ and a worst
absolute error of $0.000195$ (Section~\ref{sec:analytic}).  The second stage fits coefficients to the median.  Each fitted constant
buys accuracy at the cost of complexity, and the useful results form a Pareto frontier of worst error against number of fitted
constants.  After discussing how to measure
accuracy and complexity (Section~\ref{sec:measures}), we fit rational functions, with up to
six fitted constants, to minimize the worst absolute-or-relative error over all $\alpha > 0$.  The frontier
reaches $0.000087$ with two fitted constants and $0.000041$ with five,
mostly with some coefficients still set from $g^*$; Form C is best with zero and one fitted constants, and Form B, with
rationals of degree 4, from two to five.  The end function from the beta median is what makes the difference:  with every
scale set, Forms B and C have worst errors 6 and 12 times smaller than Form A's ($0.0021$ and $0.0011$, against $0.013$), and with
the same number of fitted constants, from two to five, Form B's are 1.8 to 2.5 times smaller than Form A's.

Section~\ref{sec:limits} collects what is known about the median at the two ends of the range, and Section~\ref{sec:forms}
defines the three forms and their ideal interpolators.  Section~\ref{sec:analytic} introduces the interpolator family and the
approximations with every coefficient set from $g^*$.  Section~\ref{sec:measures} introduces Pareto optimization and the measures
of accuracy and complexity it needs,
Section~\ref{sec:fitting} the fitting, and Section~\ref{sec:results} the Pareto frontiers, including how they depend on the
error measure.  Section~\ref{sec:numerics} describes the numerical methods, and
Section~\ref{sec:summary} summarizes the results, with the formulas of three frontier models.

\section{The median at the two ends}
\label{sec:limits}

\paragraph{Small $\alpha$.}
The median is extremely small for small $\alpha$; for example $\nuG(0.001) \approx 5.2 \times 10^{-302}$.
Its behavior there follows from the series
$P(\alpha, x) = x^\alpha M(\alpha, \alpha + 1, -x)/\Gamma(\alpha + 1)$, where $M$ is Kummer's confluent
hypergeometric function \citep[\S8.5]{dlmf}.  Setting $P = \frac{1}{2}$ at $x = \nuG$ gives
\begin{equation}
\nuG(\alpha) = 2^{-1/\alpha}\,\Gamma(1 + \alpha)^{1/\alpha}\, M(\alpha, \alpha + 1, -\nuG)^{-1/\alpha} \, .
\label{eq:smallexact}
\end{equation}
Since $M(\alpha, \alpha + 1, -x) = 1 - \alpha x/(\alpha + 1) + O(x^2)$, the last factor is $1 + O(\nuG)$, and
$\nuG = O(2^{-1/\alpha})$ is smaller than any power of $\alpha$.  So
$\nuG(\alpha) = 2^{-1/\alpha}\,\Gamma(1 + \alpha)^{1/\alpha}\,\big(1 + O(2^{-1/\alpha})\big)$.  Dropping the last factor,
which exceeds 1, gives the lower bound $2^{-1/\alpha}\,\Gamma(1 + \alpha)^{1/\alpha}$ of \citet{lyon2021closed}; and from
$\log\Gamma(1 + \alpha) = -\gamma\alpha + \sum_{j \ge 2} (-1)^j \zeta(j)\,\alpha^j/j$,
\begin{equation}
\nuG(\alpha) = 2^{-1/\alpha} e^{-\gamma}\left(1 + \frac{\pi^2}{12}\,\alpha + O(\alpha^2)\right) \quad \text{as } \alpha \to 0 \, .
\label{eq:small}
\end{equation}
The leading term is the asymptote of \citet{berg2006chen}, and the linear term is the one used by
\citet{lyon2021closed}; \eqref{eq:smallexact} shows that both are terms of the full asymptotic expansion, to all orders in $\alpha$.

\paragraph{Large $\alpha$.}
\citet{choi1994medians} and \citet{berg2006chen} found the asymptotic series
\begin{equation}
\nuG(\alpha) = \alpha - \frac{1}{3} + \frac{8}{405\,\alpha} + \frac{184}{25515\,\alpha^2} + \frac{2248}{3444525\,\alpha^3} + O(\alpha^{-4})
\quad \text{as } \alpha \to \infty \, .
\label{eq:large}
\end{equation}
The coefficient $\frac{184}{25515} = \frac{2^3 \cdot 23}{3^6 \cdot 5 \cdot 7}$ is as given by \citet{berg2006chen};
it appears as $\frac{144}{3^6 \cdot 5 \cdot 7}$ in \citet{choi1994medians}, a misprint also noted by \citet{lyon2021closed}.

\paragraph{Special values.}  At $\alpha = 1$, $\nuG(1) = \log 2$.  At $\alpha = 2$, where the distribution function is
${1 - (1 + x)e^{-x}}$, the median is $\nuG(2) = -1 - W_{-1}\big(-1/(2e)\big) \approx 1.678347$, with $W_{-1}$ the lower real branch of
the Lambert $W$ function \citep{corless1996lambert}, and at $\alpha = \frac{1}{2}$, it is
$\nuG(\frac{1}{2}) = \big(\operatorname{erf}^{-1}(\frac{1}{2})\big)^2 \approx 0.227468$ \citep{lyon2026median}.  Of these, only
$\alpha = 1$ gives a value in elementary functions, and it is the one value used below to set coefficients by exactness.

\section{Three interpolated forms}
\label{sec:forms}

\paragraph{Three exact rewritings.}  We rewrite the exact median in these three ways,
\begin{equation}
\nuG(\alpha) = \frac{A(\alpha) + \alpha}{2^{1/\alpha}} = \frac{B(\alpha)}{2^{1/\alpha} - 1}
= \frac{C(\alpha)}{2^{1/\alpha} - 1 - c/\alpha^2} \, , \qquad c = \Big(\frac{1}{2}\log 2 - \frac{1}{3}\Big)\log 2 \approx 0.0091774 \, ,
\label{eq:ABC}
\end{equation}
which define three functions $A$, $B$, and $C$ of $\alpha$; $A$ is the function of \citet{lyon2021closed}
(see Figure~\ref{fig:ABC}, top panel).  All three of these functions tend to $e^{-\gamma} \approx 0.561459$ as $\alpha \to 0$,
where each denominator is
$2^{1/\alpha}$ to within a factor $1 + O(2^{-1/\alpha}/\alpha^2)$ and \eqref{eq:small} applies.  At large $\alpha$,
$2^{1/\alpha} = 1 + \log 2/\alpha + O(\alpha^{-2})$, so subtracting 1 in the denominators of $B$ and $C$ leaves about
$\log 2/\alpha$, which moves a factor $\alpha$ into the numerator, as the explicit $+\alpha$ does for $A$.  With
\eqref{eq:large}, $A$ falls overall to $\Ainf = \log 2 - \frac{1}{3} \approx 0.359814$, while
$B$ and $C$ rise overall to $\log 2 \approx 0.693147$, as
\[
B(\alpha) = \log 2 + \frac{c}{\alpha} + O(\alpha^{-2}) \, , \qquad
C(\alpha) = \log 2 - \Big(\frac{37}{405} - \frac16\log^2 2\Big)\frac{\log 2}{\alpha^2} + O(\alpha^{-3}) \, :
\]
the constant $c$ is the $1/\alpha$ coefficient of $B$, and the $c/\alpha^2$ in the last denominator removes it, so that $C$
approaches $\log 2$ one order faster.  What is known about bounds is a statement about these functions.  $A$ lies strictly
between its limits for all $\alpha > 0$ (\citealp{lyon2023tight}; Appendix~\ref{app:Abounds} gives a shorter proof), and numerically it decreases monotonically, the open conjecture of
Section~\ref{sec:intro}.  Since $\nu_{\Gamma 4} = \log 2/(2^{1/\alpha} - 1)$, Theorem~3 of \citet{lyon2026median} says that
$B < \log 2$ for $\alpha < 1$ and $B > \log 2$ for $\alpha > 1$; numerically $B$ rises to $0.695203$ near $\alpha = 2.14$
before falling back to $\log 2$.  $C$ lies strictly between its limits for all $\alpha > 0$
(Theorem~\ref{prop:Cbounds}, Appendix~\ref{app:Cbounds}), and numerically it is strictly increasing; we have not proved
that, and it neither implies nor follows from the conjecture for $A$ in any direct way (Remark~\ref{rem:AC}).
Neither monotonicity is within reach of the methods of Appendices~\ref{app:Abounds} and~\ref{app:Cbounds}, which bound $A$
and $C$ by comparing $\nuG$ with explicit bounds on it.  At large $\alpha$,
\[
A'(\alpha) \approx -\Big(\frac{1}{2}\log^2 2 - \frac13\log 2 + \frac{8}{405}\Big)\frac{1}{\alpha^2} \approx -\frac{0.028931}{\alpha^2} \, ,
\qquad C'(\alpha) \approx \frac{0.015641}{\alpha^3} \, ,
\]
so a proof of either would need two-sided bounds on $\nuG'$ accurate to $o(\alpha^{-2})$, far beyond the known
$\nuG' < 1$ and convexity of $\nuG$.

\paragraph{Interpolated forms.}  Each form approximates its function by interpolating between the function's two limits,
with an interpolator $g(\alpha)$ of its own that goes from 0 at $\alpha = 0$ to 1 as $\alpha \to \infty$:
\[
\hat A = (1 - g)\,e^{-\gamma} + g\,\Ainf \, , \qquad \hat B = (1 - g)\,e^{-\gamma} + g\log 2 \, , \qquad \hat C = (1 - g)\,e^{-\gamma} + g\log 2 \, .
\]
In \eqref{eq:ABC} these give the three approximations $\nuA$, $\nuB$, and $\nuC$, each a single formula, one smooth expression
for all $\alpha > 0$.  That is, each interpolation between small-$\alpha$ and large-$\alpha$ limits of the numerator functions is
equivalent to an interpolation between small-$\alpha$ and large-$\alpha$ function approximations to the median
\begin{equation}
\nuh(\alpha) = \big(1 - g(\alpha)\big)\,\nu_0(\alpha) + g(\alpha)\,\nu_\infty(\alpha) \, ,
\label{eq:interp}
\end{equation}
where the \emph{end functions} $\nu_0$ and $\nu_\infty$ are \eqref{eq:ABC} with the function replaced by its small-$\alpha$
and large-$\alpha$ limits.

Form A is the form of \citet{lyon2021closed}, and its two ends are the tightest upper and lower bounds
of the form $2^{-1/\alpha}(A + \alpha)$ with constant $A$ \citep{lyon2023tight}.  Form B's $\nu_\infty$ is $\nu_{\Gamma 4}$,
which is exact at $\alpha = 1$, a bound on each side of it, and approaches $\alpha - \frac{1}{2}\log 2$ at large $\alpha$, about
$0.0132403$ below the true $\alpha - \frac{1}{3}$.  Form B's main role here is as the step to Form C, though with enough fitted
constants it does best (Section~\ref{sec:results}).  Form C's $\nu_\infty$ is $\log 2/(2^{1/\alpha} - 1 - c/\alpha^2) =
\alpha - \frac{1}{3} + (\frac19 - \frac16\log^2 2)/\alpha + O(\alpha^{-2})$, whose $1/\alpha$ coefficient $0.031036$ differs from
the $\frac{8}{405} \approx 0.019753$ of \eqref{eq:large}; its denominator is positive for all $\alpha > 0$.  For Forms B and C,
$\nu_0$ is just $\nu_\infty$ scaled down by the factor $e^{-\gamma}/\log 2$.

\paragraph{Convex combinations.}  With $0 \le g \le 1$, each interpolated function is a convex combination of its two limits,
and \eqref{eq:interp} is a convex combination of $\nu_0$ and $\nu_\infty$.  The weight $g^*$ that makes it exact
(Section~\ref{sec:ideal}) lies between 0 and 1 exactly where the function lies between its limits, that is, where $\nu_0$ and
$\nu_\infty$ are a lower and an upper bound on the median; any interpolator with values between 0 and 1 then gives a value
between the bounds.  That holds for all $\alpha$ for Forms A and C, as Theorems~\ref{prop:Abounds} and~\ref{prop:Cbounds} show.  Form B is the exception:
$B$ exceeds its limit $\log 2$ for $\alpha > 1$, so its exact weight exceeds 1 there, and Form B's interpolator must extrapolate
beyond $\nu_\infty$ rather than average.

\begin{figure}[tbp]
\centering
\includegraphics[width=0.9\linewidth]{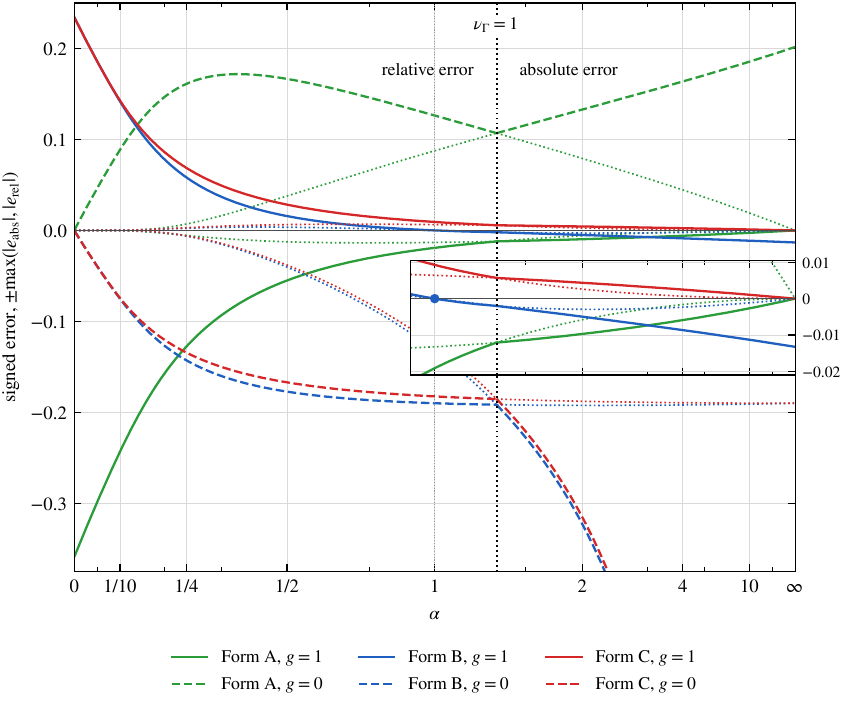}
\caption{Errors of the end functions alone, $\nu_0$ ($g = 0$, dashed) and $\nu_\infty$ ($g = 1$, solid), for each form,
signed so that positive values are upper bounds.
The solid and dashed lines show the measure \eqref{eq:loss} of Section~\ref{sec:measures}:  the relative error left of the black
dotted line at $\alpha = 1.31425$, where $\nuG = 1$, and the absolute error right of it.  Thin dotted lines continue the other error,
where it is the smaller.  The inset magnifies the errors for $\alpha > 0.9$ by $4\times$ vertically, on the same $\alpha$ scale.
Form A's ends are the bounds of \citet{lyon2023tight}; Form B's $\nu_\infty$ is $\nu_{\Gamma 4}$
of \citet{lyon2026median}, which crosses the true median at $\alpha = 1$ (grey dotted line; blue dot in the inset).  With $g = 0$,
Forms B and C have slope $e^{-\gamma}/\log 2$ at large $\alpha$, so their errors are unbounded as $\alpha \to \infty$;
with $g = 1$, the errors of Forms A and C go to 0 there, while Form B's absolute error approaches $\frac{1}{3} - \frac{1}{2}\log 2$.
}
\label{fig:ends}
\end{figure}

\paragraph{Absolute and relative errors.}
We write $e_{\mathrm{abs}} = \nuh - \nuG$ and $e_{\mathrm{rel}} = \nuh/\nuG - 1$ for the absolute and relative errors, which
always have the same sign.  In all three forms, $\nu_0 = 2^{-1/\alpha}e^{-\gamma}(1 + O(\alpha))$ as $\alpha \to 0$, so the
relative error vanishes there if $g(0) = 0$; and at large $\alpha$ the absolute error vanishes if $g \to 1$ along an appropriate trajectory.
Figure~\ref{fig:ends} shows the errors of the end functions alone.

\subsection{Ideal interpolators}
\label{sec:ideal}

For each form, the ideal interpolator $g^*(\alpha) = (\nuG - \nu_0)/(\nu_\infty - \nu_0)$ makes \eqref{eq:interp} exact
(Figure~\ref{fig:ABC}, bottom); it is the function of \eqref{eq:ABC} rescaled, $g^* = (A - e^{-\gamma})/(\Ainf - e^{-\gamma})$ for
Form A and $(B - e^{-\gamma})/(\log 2 - e^{-\gamma})$ and $(C - e^{-\gamma})/(\log 2 - e^{-\gamma})$ for Forms B and C.  From \eqref{eq:small} and \eqref{eq:large}, its behavior at the ends is:
\begin{itemize}
\item Form A\@:  $g^* = 2.669128\,\alpha + O(\alpha^2)$ as $\alpha \to 0$, and $1 - g^* \sim 0.143472/\alpha$ as
$\alpha \to \infty$, with $g^*(1) = (1 + e^{-\gamma} - 2\log 2)/(e^{-\gamma} - \Ainf) = 0.868678$ (the quantities $P_0$,
$P_\infty$, and $P_1$ of \citet{lyon2021closed}), where $2.669128 = (1 - \frac{\pi^2}{12}e^{-\gamma})/(e^{-\gamma} - \Ainf)$ and
$0.143472 = (\frac{8}{405} + \frac{1}{2}\log^2 2 - \frac{1}{3}\log 2)/(e^{-\gamma} - \Ainf)$.
\item Forms B and C\@:  $g^* = 3.506644\,\alpha + O(\alpha^2)$ as $\alpha \to 0$, where $3.506644 = \pi^2 e^{-\gamma}/(12(\log 2 - e^{-\gamma}))$.
\item Form B\@:  $g^*(1) = 1$ exactly, and $g^* - 1 \sim 0.069691/\alpha$ as $\alpha \to \infty$, where
$0.069691 = (\frac{1}{2}\log 2 - \frac{1}{3})\log 2/(\log 2 - e^{-\gamma})$:  the ideal interpolator exceeds 1 for $\alpha > 1$ and
approaches it from above.
\item Form C\@:  $1 - g^* \sim 0.059386/\alpha^2$ as $\alpha \to \infty$, where
$0.059386 = (\frac{37}{405} - \frac16\log^2 2)\log 2/(\log 2 - e^{-\gamma})$, with $\frac{37}{405} = \frac19 - \frac{8}{405}$
from the two $1/\alpha$ coefficients above; and $g^*(1) = ((1 - c)\log 2 - e^{-\gamma})/(\log 2 - e^{-\gamma}) = 0.951694$.
\end{itemize}
Form B's ideal interpolator thus approaches 1 from above, as its ends require:  an interpolator that never exceeds 1, such as
an arctangent, cannot make Form B's absolute error vanish at large $\alpha$, and its worst error is at least $\frac{1}{2}\log 2 - \frac{1}{3} \approx 0.0132403$.

\section{Analytic interpolators}
\label{sec:analytic}

Before fitting anything, we ask how well an interpolator can do with only what is known about $g^*$ in closed form.

\paragraph{Interpolator family.}  All the interpolators we report are rational functions built from polynomials with positive
coefficients, which we write in what we call division form:  the coefficient of $\alpha^k$ is written $a^{-k}$, so that the
term is $(\alpha/a)^k$ and its \emph{scale} $a$ is the value of $\alpha$ at which the term reaches 1.  Scales make fits better
conditioned than coefficients, and provide a more intuitive numerical interpretation.  Positive coefficients keep every
interpolator positive, and every denominator free of zeros,
for all $\alpha > 0$, with nothing to check.  The numerator is a polynomial
\[
S(\alpha) = \sum_{i \in N} \left(\frac{\alpha}{a_i}\right)^i ,
\]
whose powers $N$ are its degree $m$ and any of the lower powers; with all of them, $S = \alpha/a_1 + (\alpha/a_2)^2 + \cdots +
(\alpha/a_m)^m$.  The interpolator is
\[
g(\alpha) = \frac{S(\alpha)}{1 + \sum_{i \in D} (\alpha/b_i)^i + (\alpha/a_m)^m} \, ,
\]
whose denominator shares the top term of $S$, so that $g(\infty) = 1$, and has each lower power $i < m$ absent, free (its own
scale $b_i$), or, if $S$ has that power, tied ($b_i = a_i$, the same term as in $S$).  With every power in $S$ and every lower
term tied, it is $S/(1 + S)$.  A rational can exceed 1 and approach 1 from above at large $\alpha$, as Form B requires.  The
numerators with some lower power left out matter:  the best rationals with two or more fitted constants have no $\alpha^2$ term
in the numerator (Section~\ref{sec:results}).  Letting the middle coefficients of the numerators with every power take either
sign, with the denominator kept positive, does not lower the best worst errors of Section~\ref{sec:results} (see the notebook).

\paragraph{Arctangents.}  For Form A, \citet{lyon2021closed} also used the one-scale arctangent interpolator
$g = \frac{2}{\pi}\arctan(\alpha/a)$, whose scale $a$ is the pivot where $g = \frac{1}{2}$.  It does well for a single
constant:  with $a = 0.209257$, set by exactness at $\alpha = 1$, its worst error \eqref{eq:loss} is $0.00613$;
and with $a$ fitted, $0.00396$, at $a = 0.216391$.
These were better than the best one-scale rational-function interpolators compared for Form A \citep{lyon2021closed}.
Here we also fitted arctangents $\frac{2}{\pi}\arctan S$ of the polynomials above, of degree up to 4, and arctangents
of rational functions, but leave them out as none came close to the best rational-function interpolators.

\paragraph{Conditions from $g^*$.}  Three properties of each form's ideal interpolator are known in closed form
(Section~\ref{sec:forms}):  its slope $g^{*\prime}(0)$, the coefficient $K$ in $1 - g^* \sim K/\alpha$ for Forms A and B
($K < 0$ for Form B) or $K/\alpha^2$ for Form C, and its value $g^*(1)$.  An interpolator can take any of them in place of a
scale fitted to the median.  For a rational $g = P/Q$, with coefficients $u = a_i^{-i}$ or $b_i^{-i}$, each is a linear
condition on the coefficients:  $g'(0) = u_{a_1}$; the large-$\alpha$ coefficient of $1 - g = (Q - P)/Q$ is the coefficient of
the right power of $\alpha$ in $Q - P$, divided by $u_{a_m}$; and $g(1) = g^*(1)$ is $P(1) = g^*(1)\,Q(1)$.  Each condition is
solved for one coefficient, which must come out positive, and sets one scale; the first two make the error vanish one order
faster at that end.  These are the conditions of \citet{lyon2021closed}, which used all three for Form A with
rationals in coefficient form.  Its single-scale rationals ``best at low $k$'' and ``exact at $k = 1$'' ($b_0 = 0.374654$ and
$0.151175$) are ours with the slope set and with exactness at $\alpha = 1$, and its degree-3 rational with all three set,
which has two more scales to fit, is still the best of Form A with two fitted constants (Section~\ref{sec:results}).

\paragraph{Every scale set.}  When a structure has as many scales as the conditions it takes, every scale is set:  the
interpolator is optimized analytically, by matching $g^*$ where $g^*$ is known, with no error optimization criterion, and its scales are
closed-form functions of $e^{-\gamma}$, $\log 2$, and $\pi$.  Table~\ref{tab:analytic} lists, for each form and each combination
of conditions, the structure with the smallest worst relative error, with the end functions alone for comparison.  For each it
gives the worst relative error and the worst absolute error over all $\alpha > 0$, including the limits
as $\alpha \to 0$ and $\alpha \to \infty$.

\begin{table}[tbp]
\centering
\footnotesize
\caption{Analytic interpolators.  For each form and each combination of the conditions from $g^*$ (value $g^*(1)$, slope
$g^{*\prime}(0)$, and large-$\alpha$ coefficient $K$), the rational interpolator with every scale set from them
and the smallest worst relative error, after the end functions $g = 0$ and $g = 1$.  Worst errors are over all $\alpha > 0$,
with the $\alpha$ at which they occur ($0$ and $\infty$ are limits).  The combinations not listed admit no structure in our family with
every scale positive.}
\label{tab:analytic}
\setlength{\tabcolsep}{2pt}
\renewcommand{\arraystretch}{1.35}
\begin{tabular}{@{}l >{\raggedright}p{0.08\linewidth} l >{\raggedright}p{0.27\linewidth}@{\hspace{6pt}}l l@{}}
\toprule
Form & Set from & Interpolator $g(\alpha)$ & Scales & \multicolumn{1}{c}{\shortstack{worst $|e_{\mathrm{rel}}|$\\(at $\alpha$)}}
  & \multicolumn{1}{c}{\shortstack{worst $|e_{\mathrm{abs}}|$\\(at $\alpha$)}} \\
\midrule
Form A & -- & $g = 0$ & -- & $0.172$ \ ($0.37$) & $0.202$ \ ($\infty$) \\
 & -- & $g = 1$ & -- & $0.359$ \ ($0$) & $0.0137$ \ ($0.75$) \\
 & $g^*(1)$ & $\dfrac{\alpha/a_1}{1 + \alpha/a_1}$ & $a_{1} = 0.151175$ & $0.0444$ \ ($0.083$) & $0.00101$ \ ($0.35$) \\
 & $g^{*\prime}(0)$ & $\dfrac{\alpha/a_1}{1 + \alpha/a_1}$ & $a_{1} = 0.374654$ & $0.0471$ \ ($0.30$) & $0.0144$ \ ($1.2$) \\
 & $K$ & $\dfrac{\alpha/a_1}{1 + \alpha/a_1}$ & $a_{1} = 0.143472$ & $0.0485$ \ ($0.084$) & $0.00135$ \ ($0.40$) \\
 & $g^*(1)$, $g^{*\prime}(0)$ & $\dfrac{\alpha/a_1 + (\alpha/a_2)^2}{1 + \alpha/a_1 + (\alpha/a_2)^2}$ & $a_{1} = 0.374654$, $a_{2} = 0.503426$ & $0.0235$ \ ($0.24$) & $0.00378$ \ ($3.1$) \\
 & $K$, $g^*(1)$ & $\dfrac{(\alpha/a_2)^2}{1 + \alpha/b_1 + (\alpha/a_2)^2}$ & $a_{2} = 0.0877635$, $b_{1} = 0.0536859$ & $0.0160$ \ ($0.12$) & $0.000474$ \ ($0.34$) \\
 & $K$, $g^*(1)$, $g^{*\prime}(0)$ & $\dfrac{\alpha/a_1 + (\alpha/a_3)^3}{1 + \alpha/a_1 + (\alpha/b_2)^2 + (\alpha/a_3)^3}$ & $a_{1} = 0.374654$, $a_{3} = 0.234607$,\newline $b_{2} = 0.300004$ & $0.0131$ \ ($0.16$) & $0.000527$ \ ($0.34$) \\
\midrule
Form B & -- & $g = 0$ & -- & $0.192$ \ ($2.1$) & $\infty$ \\
 & -- & $g = 1$ & -- & $0.235$ \ ($0$) & $0.0132$ \ ($\infty$) \\
 & $g^{*\prime}(0)$ & $\dfrac{\alpha/a_1}{1 + \alpha/a_1}$ & $a_{1} = 0.285173$ & $0.0552$ \ ($0.41$) & $0.0674$ \ ($\infty$) \\
 & $K$, $g^*(1)$ & $\dfrac{\alpha/a_1 + (\alpha/a_2)^2}{1 + (\alpha/a_2)^2}$ & $a_{1} = 1$, $a_{2} = 0.263991$ & $0.0301$ \ ($0.13$) & $0.00139$ \ ($0.38$) \\
 & $K$, $g^{*\prime}(0)$ & $\dfrac{\alpha/a_1 + (\alpha/a_2)^2}{1 + (\alpha/a_2)^2}$ & $a_{1} = 0.285173$, $a_{2} = 0.140975$ & $0.0590$ \ ($0.17$) & $0.00672$ \ ($0.75$) \\
 & $K$, $g^*(1)$, $g^{*\prime}(0)$ & $\dfrac{\alpha/a_1 + (\alpha/a_2)^2}{1 + \alpha/b_1 + (\alpha/a_2)^2}$ & $a_{1} = 0.285173$, $a_{2} = 0.263991$,\newline $b_{1} = 0.398940$ & $0.00211$ \ ($0.077$) & $0.000185$ \ ($0.44$) \\
\midrule
Form C & -- & $g = 0$ & -- & $0.190$ \ ($\infty$) & $\infty$ \\
 & -- & $g = 1$ & -- & $0.235$ \ ($0$) & $0.00687$ \ ($0.69$) \\
 & $g^*(1)$ & $\dfrac{(\alpha/a_2)^2}{1 + (\alpha/a_2)^2}$ & $a_{2} = 0.225295$ & $0.0388$ \ ($0.11$) & $0.00130$ \ ($0.36$) \\
 & $g^{*\prime}(0)$ & $\dfrac{\alpha/a_1}{1 + \alpha/a_1}$ & $a_{1} = 0.285173$ & $0.0434$ \ ($0.40$) & $0.0542$ \ ($\infty$) \\
 & $K$ & $\dfrac{(\alpha/a_2)^2}{1 + (\alpha/a_2)^2}$ & $a_{2} = 0.243693$ & $0.0443$ \ ($0.12$) & $0.00215$ \ ($0.43$) \\
 & $g^*(1)$, $g^{*\prime}(0)$ & $\dfrac{\alpha/a_1 + (\alpha/a_2)^2}{1 + \alpha/a_1 + (\alpha/a_2)^2}$ & $a_{1} = 0.285173$, $a_{2} = 0.248493$ & $0.00204$ \ ($0.19$) & $0.000112$ \ ($0.37$) \\
 & $K$, $g^*(1)$ & $\dfrac{\alpha/a_1 + (\alpha/a_2)^2}{1 + \alpha/a_1 + (\alpha/a_2)^2}$ & $a_{1} = 0.349355$, $a_{2} = 0.243693$ & $0.00696$ \ ($0.14$) & $0.000297$ \ ($0.37$) \\
 & $K$, $g^{*\prime}(0)$ & $\dfrac{\alpha/a_1 + (\alpha/a_2)^2}{1 + \alpha/a_1 + (\alpha/a_2)^2}$ & $a_{1} = 0.285173$, $a_{2} = 0.243693$ & $0.00111$ \ ($0.18$) & $0.000195$ \ ($1.1$) \\
 & $K$, $g^*(1)$, $g^{*\prime}(0)$ & $\dfrac{\alpha/a_1 + (\alpha/a_4)^4}{1 + (\alpha/b_2)^2 + (\alpha/a_4)^4}$ & $a_{1} = 0.285173$, $a_{4} = 0.238100$,\newline $b_{2} = 0.232636$ & $0.0101$ \ ($0.18$) & $0.000523$ \ ($0.36$) \\
 
\bottomrule
\end{tabular}
\end{table}

The end functions alone do poorly (Figure~\ref{fig:ends}).  In Form A they are the bounds of \citet{lyon2023tight}:  the upper
bound $2^{-1/\alpha}(e^{-\gamma} + \alpha)$ has worst absolute error $e^{-\gamma} - \Ainf \approx 0.201646$, approached as
$\alpha \to \infty$, and the lower bound $2^{-1/\alpha}(\Ainf + \alpha)$ has worst relative error
$1 - e^{\gamma}\Ainf \approx 0.359145$, approached as $\alpha \to 0$.  In Forms B and C, $g = 1$ has worst relative error
$e^{\gamma}\log 2 - 1 \approx 0.234545$ as $\alpha \to 0$, and $g = 0$ leaves the slope wrong at large $\alpha$, so its absolute
error is unbounded.

A single condition helps less than one might hope.  Each possible condition sets the scale of a one-scale rational,
$(\alpha/a_1)/(1 + \alpha/a_1)$ or $(\alpha/a_2)^2/(1 + (\alpha/a_2)^2)$, and leaves worst relative errors of $0.039$ to $0.055$.  For Form A, exactness at $\alpha = 1$
gives $a_1 = 1/g^*(1) - 1 = 0.151175$, the rational of \citet{lyon2021closed} exact at $\alpha = 1$, with $0.0444$.  For Form B,
$g^*(1) = 1$ forces $g = 1$.  The slope or the large-$\alpha$ coefficient alone corrects the error at one end only.

\begin{figure}[tbp]
\centering
\includegraphics[width=0.9\linewidth]{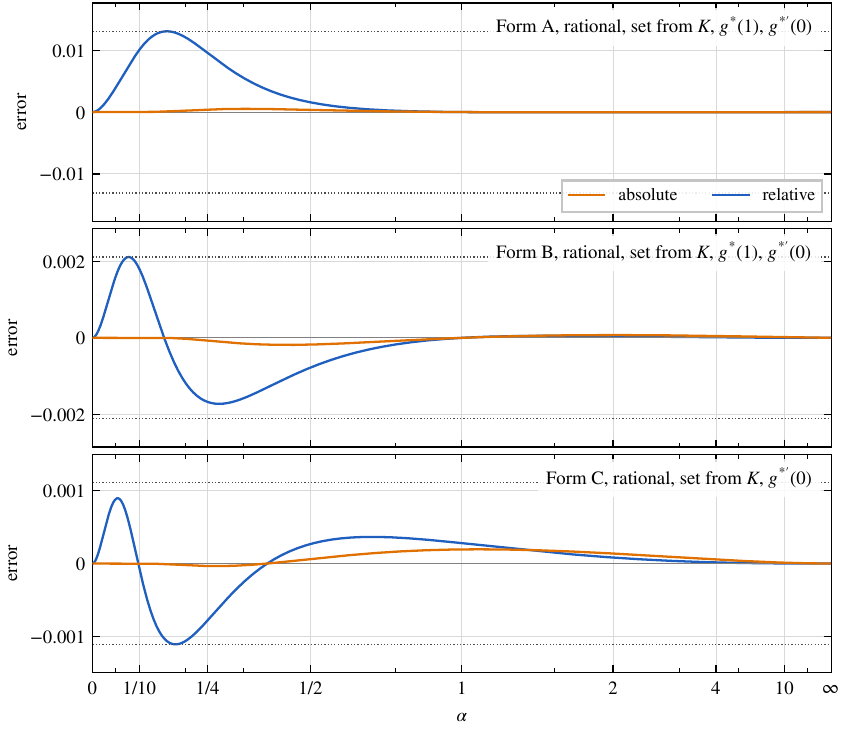}
\caption{Absolute (orange) and relative (blue) errors of the best analytic interpolator of each form
(Table~\ref{tab:analytic}), every scale set from $g^*$, with dotted lines at $\pm$ the worst relative error.}
\label{fig:analytic}
\end{figure}

With two conditions, Form C's rational $g = S/(1 + S)$, $S = \alpha/a_1 + (\alpha/a_2)^2$, has $g'(0) = 1/a_1$ and
$1 - g \sim a_2^2/\alpha^2$, so the slope $g^{*\prime}(0) = 3.506644$ and the coefficient $K = 0.059386$ of
Section~\ref{sec:forms} give $a_1 = 1/g^{*\prime}(0) = 0.285173$ and $a_2 = K^{1/2} = 0.243693$.  Its worst relative error,
$0.00111$, is the smallest of all the analytic interpolators, and its worst absolute error is $0.000195$
(Figure~\ref{fig:analytic}).
That rational is already nearly exact at $\alpha = 1$ ($g(1) = 0.95315$, against $g^*(1) = 0.951694$),
and the best rational that takes exactness there as well, $(\alpha/a_1 + (\alpha/a_4)^4)/(1 + (\alpha/b_2)^2 + (\alpha/a_4)^4)$,
has degree 4 and gives a worse $0.0101$.
For Form B, a degree-2 rational with a free denominator term takes all three conditions and gives a worst error of
$0.00211$.  Form A does least well:  its best, a degree-3 rational with no $\alpha^2$ term in its numerator and all three
conditions set, gives $0.0131$.

With any condition besides the slope, the worst error is relative, at $\alpha \le 0.24$, where the median is below $0.04$,
and it exceeds the worst absolute error by a factor of $5.7$ to $44$.  For these interpolators, the choice of whether and how to
include the absolute error hardly matters.

\section{Pareto optimization}
\label{sec:measures}

The analytic interpolators discussed above use only what is known about the median in closed form.
To do better in terms of accuracy, we fit some scales to the median numerically:  each fitted constant can buy accuracy,
at the cost of a more complex formula.
The idea of Pareto optimization is to find those approximations that cannot be beaten in both respects, accuracy and complexity,
by any other.  The Pareto frontier is the collection of these:  for each level of complexity, the one with the best accuracy.
Drawing the frontier needs specific definitions for both accuracy and complexity, the two criteria to optimize and trade.

\paragraph{Accuracy.}  The absolute and relative errors suit different ends of the range.  At small $\alpha$ the median is
tiny, and only the relative error is meaningful; at large $\alpha$ the relative error of any of the forms tends to zero, and the
absolute error is what shows whether the offset $-\frac{1}{3}$ is right.  The worst relative error alone ignores the large-$\alpha$
end, and the worst absolute error alone ignores the small-$\alpha$ end.  Their pointwise maximum weighs the two ends, and that is
what we have chosen to minimize:
\begin{equation}
\max_{\alpha > 0}\, \max\big(|e_{\mathrm{abs}}|, |e_{\mathrm{rel}}|\big) = \max_{\alpha > 0}\, \big(|e_{\mathrm{rel}}|\,\max(1, \nuG)\big)
= \max_{\alpha > 0}\, \big(|e_{\mathrm{abs}}| \big/ \min(1, \nuG)\big) \, ,
\label{eq:loss}
\end{equation}
which is the relative error where the median is below 1 ($\alpha < 1.31425$) and the absolute error above.  Their pointwise
sum, within a factor of 2 of the maximum, would serve as well.  A single number that bounds both errors is a common specification
for numerical approximations.  For
the analytic interpolators, \eqref{eq:loss} is the worst relative error except for Forms B and C where only the slope at 0 is set
(Table~\ref{tab:analytic}).  Numerically fitted interpolators balance their errors across the range, and there
the measures differ more; Section~\ref{sec:otherloss} compares the frontiers they give.

\paragraph{Complexity.}  The complexity of an approximation, to be traded against accuracy for a Pareto optimization, can
be measured in several ways.  We count fitted constants $n$, the scales fitted to the median.
Tied terms and scales set from $g^*$ do not count by this measure, so the analytic interpolators have $n = 0$.
This measure counts how much must be learned from the median beyond its exact value at $\alpha = 1$ and its asymptotic expansions,
and $n$ is also the dimension of the search space in the fit.  Counting every scale instead
measures the size of the formula, and Figure~\ref{fig:paretoset} also shows the frontier in those terms.  Other possible
complexity measures are not considered in this study, though computational cost goes roughly with the count of all
scales, not with $n$ (Section~\ref{sec:summary}).

\section{Fitting}
\label{sec:fitting}

\paragraph{Structures.}  We used every rational of Section~\ref{sec:analytic} of degree 1 to 4 with at most six fitted
constants:  1, 5, 25, and 124 structures of degrees 1, 2, 3, and 4, 155 in all.  (A structure of degree 4 with every power in
its numerator and all three lower denominator terms free would have seven.)  Every structure was fitted for
each form with all its scales fitted:  471 fits, including the six end functions.  Every structure was also fitted with every
combination of the conditions of Section~\ref{sec:analytic} that it admits:  1377 fits (514, 461, and 402 for Forms A, B, and C),
35 of them with every scale set.  A combination is
admitted when its conditions can be solved for distinct coefficients:  the slope needs the numerator term of degree 1;
the large-$\alpha$ coefficient needs a term of degree $m - 1$ in $Q - P$ for Forms A and B, and for Form C
none of degree $m - 1$ (that power tied, or in neither $P$ nor $Q$) and one of degree $m - 2$; and exactness at $\alpha = 1$ needs one more
coefficient for which the linear system stays nonsingular.  During a fit, parameter values for which a solved coefficient is not
positive are rejected.  We call interpolators with any scale set from $g^*$ ``set'' and the rest ``all fitted''.

\paragraph{Method.}  Each structure was fitted by global search followed by local polishing on a logarithmic grid of $\alpha$
from $10^{-3}$ to $10^{30}$, then checked on a much denser grid, below $\alpha = 10^{-3}$ and far beyond $10^{30}$
(Section~\ref{sec:numerics}).  A minimax fit with $n$ effective fitted constants typically equioscillates, with $n + 1$
alternating extrema of equal magnitude \citep{trefethen2019approximation}.  The frontier fits reach their worst error at $n + 1$
or more values of $\alpha$, though for some of those with scales set from $g^*$, two neighboring extrema have the same sign
(Figure~\ref{fig:profiles}).

\section{Results}
\label{sec:results}

\begin{figure}[tbp]
\centering
\includegraphics[width=0.9\linewidth]{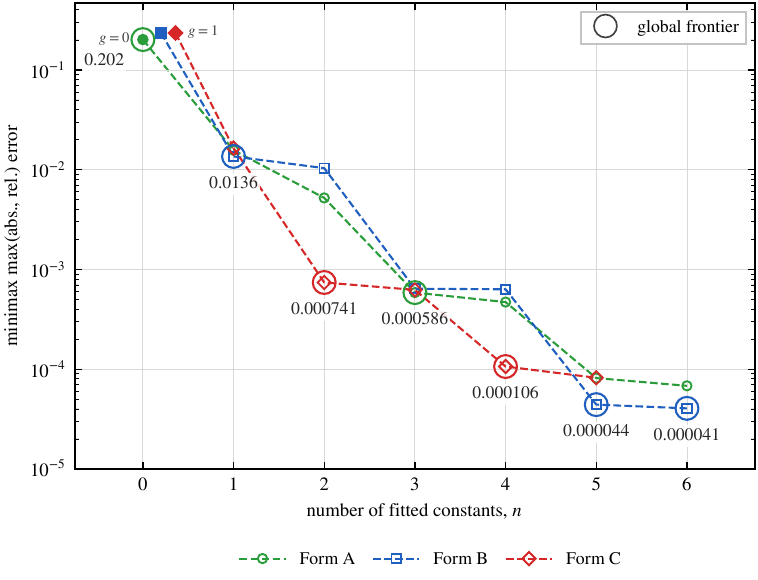}
\caption{Pareto frontier with all scales fitted:  worst error \eqref{eq:loss} against number of fitted constants, for the best
structure of each form.  Circles, in the color of the form there, mark the global
frontier, labeled with its worst errors.  At $n = 0$, each form's better end function is shown:  in each case it is the higher one,
$g = 0$ for Form A, which starts the frontier, and $g = 1$ for Forms B and C, and it is an upper bound for Forms A and C.  Each line stops at the last number of fitted constants at which its form
improves by more than 1\%; further fitted constants go unused.}
\label{fig:pareto}
\end{figure}

\begin{figure}[tbp]
\centering
\includegraphics[width=0.9\linewidth]{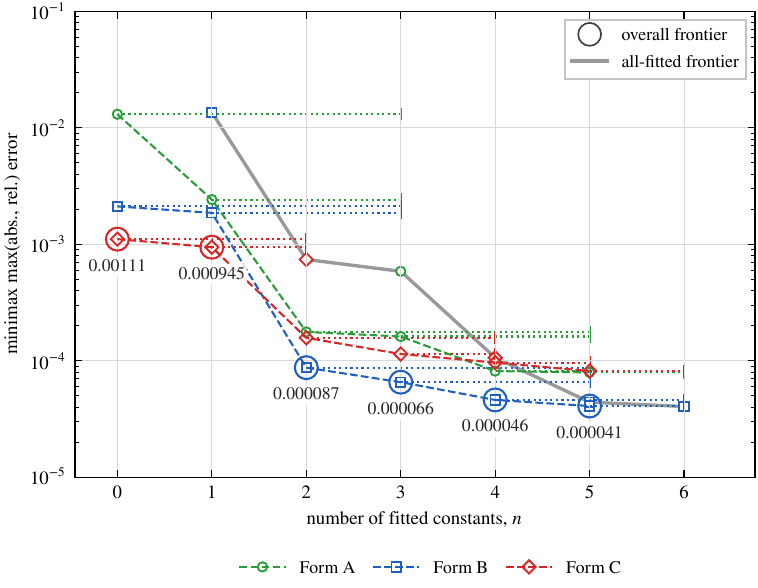}
\caption{Pareto frontier with scales set from $g^*$ (Section~\ref{sec:analytic}):  as Figure~\ref{fig:pareto}, for the best
structure of each form with some of its scales set, with the all-fitted
frontier of Figure~\ref{fig:pareto} in gray, its points in their forms' symbols.  A dotted line from each point runs to the number of constants it would have if its
set scales were counted as fitted; each ends at or above the gray frontier there, as it must.  Circles, in the color of the form there, mark the
overall frontier, of both kinds of interpolator, labeled with its worst errors.  The points with no fitted constants are the analytic
interpolators of Section~\ref{sec:analytic}.}
\label{fig:paretoset}
\end{figure}

\begin{figure}[p]
\centering
\includegraphics[width=0.86\linewidth]{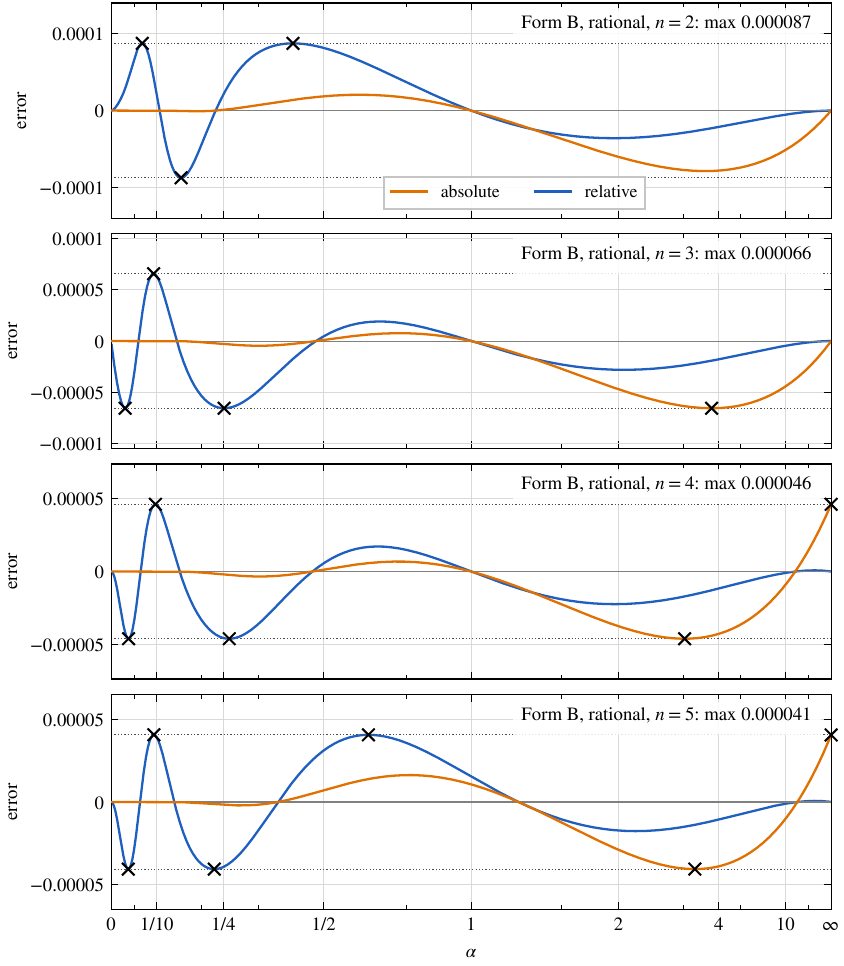}
\caption{Error profiles of the frontier models with two to five fitted constants, from Table~\ref{tab:frontier}:  the absolute
error (orange) and relative error (blue), with crosses at the worst points of the larger of the two, the error \eqref{eq:loss}
that was minimized, and dotted lines at $\pm$ its worst value; a cross at $\alpha = \infty$ marks a worst error approached in
that limit.
The profiles are for the fitted constants, before rounding to the values printed
in Table~\ref{tab:frontier}.  Each fit reaches its worst error at $n + 1$ or more values of $\alpha$.}
\label{fig:profiles}
\end{figure}

Figure~\ref{fig:pareto} gives the best worst error for each form and number of fitted constants with all scales
fitted, and Figure~\ref{fig:paretoset} with scales set from $g^*$; Table~\ref{tab:pareto} gives both, and
Table~\ref{tab:frontier} the models on the overall frontier.  A model is on a frontier only if it improves on every model
with fewer fitted constants by more than 1\%.  Their error profiles are in Figure~\ref{fig:profiles}.

\begin{table}[tbp]
\centering
\small
\caption{Best worst error \eqref{eq:loss} for each form and number of fitted constants $n$; the global
frontier in bold.  Starred entries have some scales set from $g^*$, and those with no fitted constants have all of them set; and the top two rows give the end functions alone, $g = 0$ and $g = 1$.  Each column stops at the last $n$ at which its form improves by more than 1\%.}
\label{tab:pareto}
\setlength{\tabcolsep}{7pt}
\begin{tabular}{c lll}
\toprule
$n$ & Form A & Form B & Form C \\
\midrule
$g = 0$ & \mono{0.202} & $\infty$ & $\infty$ \\
$g = 1$ & \mono{0.359} & \mono{0.235} & \mono{0.235} \\
\midrule
0 & \mono{0.0131}\rlap{$^*$} & \mono{0.00211}\rlap{$^*$} & \textbf{\mono{0.00111}}\rlap{$^*$} \\
1 & \mono{0.00242}\rlap{$^*$} & \mono{0.00186}\rlap{$^*$} & \textbf{\mono{0.000945}}\rlap{$^*$} \\
2 & \mono{0.000177}\rlap{$^*$} & \textbf{\mono{0.000087}}\rlap{$^*$} & \mono{0.000158}\rlap{$^*$} \\
3 & \mono{0.000162}\rlap{$^*$} & \textbf{\mono{0.000066}}\rlap{$^*$} & \mono{0.000114}\rlap{$^*$} \\
4 & \mono{0.000082}\rlap{$^*$} & \textbf{\mono{0.000046}}\rlap{$^*$} & \mono{0.000097}\rlap{$^*$} \\
5 & \mono{0.000080}\rlap{$^*$} & \textbf{\mono{0.000041}}\rlap{$^*$} & \mono{0.000082}\rlap{$^*$} \\
6 & \mono{0.000068} &  &  \\
 
\bottomrule
\end{tabular}
\end{table}

\begin{table}[tbp]
\centering
\small
\caption{The global Pareto frontier.  Scales are given to six significant digits, or seven where six would change the
printed worst error; the worst error is for these printed values.  Starred scales are set from $g^*$, not fitted.}
\label{tab:frontier}
\renewcommand{\arraystretch}{1.6}
\setlength{\tabcolsep}{4pt}
\begin{tabular}{c >{\raggedright}p{0.09\linewidth} >{\raggedright}p{0.38\linewidth} >{\raggedright}p{0.32\linewidth} c}
\toprule
$n$ & Form & Interpolator $g(\alpha)$ & Scales & Worst error \\
\midrule
0 & Form C & $\dfrac{\alpha/a_1 + (\alpha/a_2)^2}{1 + \alpha/a_1 + (\alpha/a_2)^2}$ & $a_{1} = 0.285173$$^*$, $a_{2} = 0.243693$$^*$ & $0.00111$ \\
\addlinespace[4pt]
1 & Form C & $\dfrac{\alpha/a_1 + (\alpha/a_2)^2}{1 + \alpha/a_1 + (\alpha/a_2)^2}$ & $a_{1} = 0.285173$$^*$, $a_{2} = 0.242863$ & $0.000945$ \\
\addlinespace[4pt]
2 & Form B & $\dfrac{\alpha/a_1 + (\alpha/a_3)^3 + (\alpha/a_4)^4}{1 + (\alpha/b_2)^2 + (\alpha/b_3)^3 + (\alpha/a_4)^4}$ & $a_{1} = 0.2851729$$^*$, $a_{3} = 0.1940768$,\newline $a_{4} = 0.2054328$, $b_{2} = 0.1549773$$^*$,\newline $b_{3} = 0.2171439$$^*$ & $0.000087$ \\
\addlinespace[4pt]
3 & Form B & $\dfrac{\alpha/a_1 + (\alpha/a_3)^3 + (\alpha/a_4)^4}{1 + (\alpha/b_2)^2 + (\alpha/b_3)^3 + (\alpha/a_4)^4}$ & $a_{1} = 0.286518$, $a_{3} = 0.183877$,\newline $a_{4} = 0.197865$, $b_{2} = 0.144401$$^*$,\newline $b_{3} = 0.205410$$^*$ & $0.000066$ \\
\addlinespace[4pt]
4 & Form B & $\dfrac{\alpha/a_1 + (\alpha/a_3)^3 + (\alpha/a_4)^4}{1 + \alpha/b_1 + (\alpha/b_2)^2 + (\alpha/b_3)^3 + (\alpha/a_4)^4}$ & $a_{1} = 0.285173$$^*$, $a_{3} = 0.180045$,\newline $a_{4} = 0.196494$, $b_{1} = 5.90709$,\newline $b_{2} = 0.142495$, $b_{3} = 0.200307$$^*$ & $0.000046$ \\
\addlinespace[4pt]
5 & Form B & $\dfrac{\alpha/a_1 + (\alpha/a_3)^3 + (\alpha/a_4)^4}{1 + \alpha/b_1 + (\alpha/b_2)^2 + (\alpha/b_3)^3 + (\alpha/a_4)^4}$ & $a_{1} = 0.285173$$^*$, $a_{3} = 0.180493$,\newline $a_{4} = 0.196530$, $b_{1} = 6.40642$,\newline $b_{2} = 0.142651$, $b_{3} = 0.200966$ & $0.000041$ \\
 
\bottomrule
\end{tabular}
\end{table}

\begin{figure}[tbp]
\centering
\includegraphics[width=0.86\linewidth]{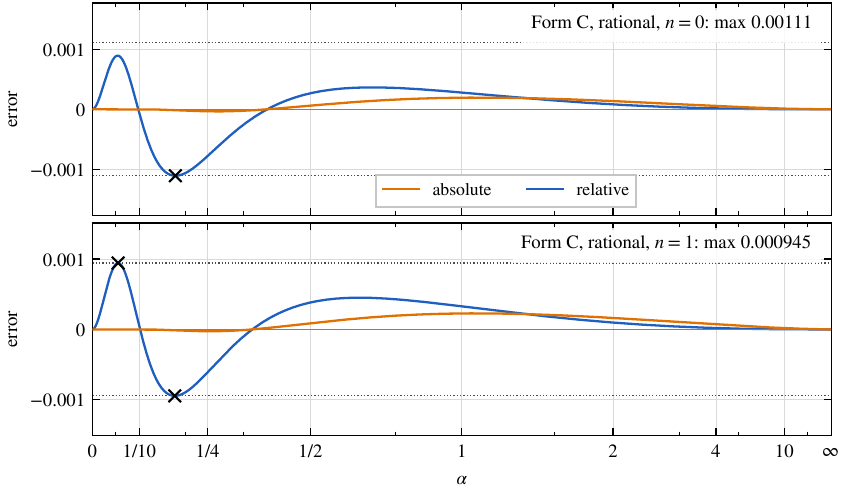}
\caption{Error profiles of the frontier models with zero and one fitted constants, Form C with $g = S/(1 + S)$,
$S = \alpha/a_1 + (\alpha/a_2)^2$, as in Figure~\ref{fig:profiles}.  With both scales set from $g^*$ (top), the worst error is a
single extremum near $\alpha = 0.18$; fitting $a_2$ (bottom) balances it against the next extremum, for a small improvement.}
\label{fig:profiles01}
\end{figure}

The frontier passes through Forms C and B, and interpolators with scales set from $g^*$ take most of it:
\begin{itemize}
\item With no fitted constants, the best analytic interpolator, Form C with $g = S/(1 + S)$, $S = \alpha/a_1 + (\alpha/a_2)^2$,
both scales set from the ends of $g^*$ (Section~\ref{sec:analytic}).
\item With one, the same with $a_1$ set and $a_2$ fitted, which balances two extrema of the error (Figure~\ref{fig:profiles01}).
With all scales fitted, the best with one is Form B's $(\alpha/a_2)^2/(1 + (\alpha/a_2)^2)$ ($0.014$).
\item With two to five, Form B rationals of degree 4,
\[
g = \frac{\alpha/a_1 + (\alpha/a_3)^3 + (\alpha/a_4)^4}{1 + \alpha/b_1 + (\alpha/b_2)^2 + (\alpha/b_3)^3 + (\alpha/a_4)^4} \, ,
\]
with no $\alpha^2$ term in the numerator, and with the term $\alpha/b_1$ only with four and five.  With two fitted constants,
all three conditions are set ($0.000087$); with three, the large-$\alpha$ coefficient $K$ and exactness at $\alpha = 1$
($0.000066$); with four, exactness at $\alpha = 1$ and the slope ($0.000046$); and with five, the slope alone ($0.000041$).
Fitting all six scales improves on five by less than 1\%.  These Form B rationals exceed 1 at large $\alpha$, as Form B
requires.
\item With all scales fitted, the best with two is Form C's $S/(1 + S)$ ($0.00074$), far better than any Form A or Form B
interpolator; with three, Form A's degree-3 rational with no $\alpha^2$ term in its numerator ($0.00059$); and with four, a
Form C degree-4 rational like the Form B ones above, with its $\alpha^3$ term tied ($0.00011$); and with five, the Form B
rational above without $\alpha/b_1$ ($0.000044$).
\item Form A is best with two fitted constants with the degree-3 rational of \citet{lyon2021closed}, with all three conditions
set; its worst error, $0.00018$, is its worst relative error, as reported there.
\end{itemize}
In every case the error at $\alpha < 10^{-3}$ is smaller than at $10^{-3}$ and tends to zero as $\alpha \to 0$, so the worst errors hold
over all $\alpha > 0$.

\begin{figure}[tbp]
\centering
\includegraphics[width=0.9\linewidth]{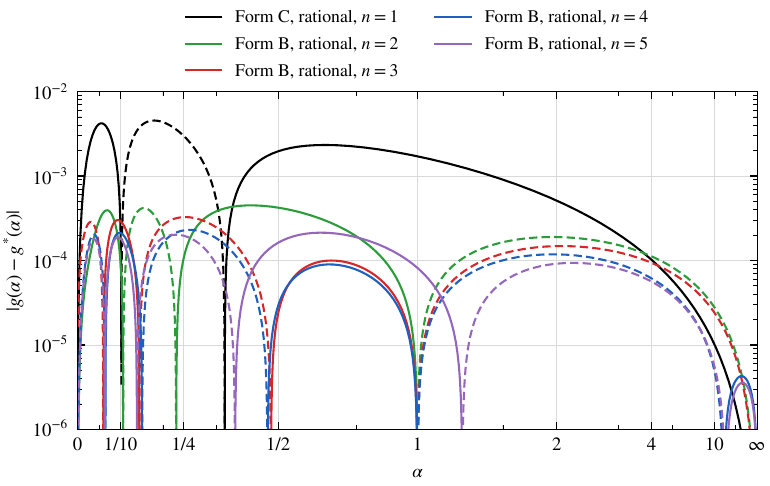}
\caption{Difference between each frontier interpolator of Table~\ref{tab:frontier} and its form's ideal interpolator $g^*$,
solid where the fit is above the ideal and dashed where below.}
\label{fig:gerror}
\end{figure}

Figure~\ref{fig:gerror} shows how closely the frontier interpolators follow the ideal:  the differences $g - g^*$ range from about
$0.0045$ for one fitted constant to about $0.0002$ for five, about five times the worst errors in the median.  The median is less
sensitive than $g$:  as $\alpha \to 0$, in Forms B and C, an error $\delta$ in $g$ changes $\nuh$ by a relative
$(e^{\gamma}\log 2 - 1)\,\delta \approx 0.23\,\delta$.

\subsection{Other error measures}
\label{sec:otherloss}

To see how the frontier depends on the choice \eqref{eq:loss}, we repeated the whole search, with all scales fitted, minimizing instead the
sum, the relative error alone, and the absolute error alone.  The sum and the relative error alone select the same frontier
models as \eqref{eq:loss} at every $n$ from 1 to 6 except 3, where both pick the Form C rational of degree 2 with $\alpha/b_1$ in its
denominator rather than the Form A rational of degree 3.  The relative error alone leaves larger absolute errors at large
$\alpha$ ($0.0018$ at $n = 3$, against $0.00049$ for the fit to \eqref{eq:loss}).  The absolute error alone gives absolute errors down to $4 \times 10^{-7}$ at $n = 6$, but ignores small
$\alpha$, where its fits have relative errors of $0.2\%$ to $4\%$.

\section{Numerical methods and use of AI tools}
\label{sec:numerics}

All computations were in double precision with NumPy \citep{harris2020array} and SciPy \citep{virtanen2020scipy}.  Medians up to
$\alpha = 10^5$ came from SciPy's inverse incomplete gamma function, which agrees with 40-digit mpmath \citep{mpmath} roots to a
relative $6 \times 10^{-13}$ or better, and above $10^5$ from the series \eqref{eq:large}.  At large $\alpha$, $\nuh - \nuG$ is a small
difference of large numbers, so it was computed from expressions rearranged in terms of $x = 1/\alpha$ and $h = 1 - g$, itself
computed without cancellation as $(Q - P)/Q$.  For Form B, for example, with $u = x\log 2$ and
$\Psi(u) = (e^u - 1 - u)/u^2$,
\[
\nuh - \nuG = -\frac{(1 - e^{-\gamma}/\log 2)\,h\,\alpha + \Psi(u)\log 2}{1 + u\,\Psi(u)} + \frac{1}{3} - T(x) \, ,
\]
where $T(x)$ is the sum of the terms of \eqref{eq:large} beyond $\alpha - \frac{1}{3}$.  For Form C,
\[
\nuh - \nuG = -\frac{(1 - e^{-\gamma}/\log 2)\,h\,\alpha + \Psi(u)\log 2 - c/\log 2}{1 + u\,\Psi(u) - (c/\log 2)\,x} + \frac{1}{3} - T(x) \, ,
\]
and for Form A,
\[
\nuh - \nuG = -\big(1 - u\,\Psi(-u)\big)\log 2 + e^{-u}\big(\Ainf + (e^{-\gamma} - \Ainf)\,h\big) + \frac{1}{3} - T(x) \, .
\]
This mattered:  some rational structures place a transition far beyond the range of interest, such as $\alpha = a_3^3/a_2^2$, which can be
$10^{20}$ with all scales between $10^{-4}$ and $10^4$.  The fitting grid therefore extends to $10^{30}$, and every fit was checked at
$10^{40}$, $10^{60}$, and $10^{100}$.  Below $\alpha = 10^{-3}$, relative errors came from \eqref{eq:smallexact} without its last factor.

Each structure was fitted by differential evolution \citep{storn1997differential} from three or more random starts, on 600
logarithmically spaced values of $\alpha$ from $10^{-3}$ to $10^5$ and 60 from $10^5$ to $10^{30}$, then polished by Nelder--Mead
\citep{nelder1965simplex}, with refinement passes adding the worst points of a check grid of about $10^5$ values.  The search
ranges were $10^{-4}$ to $10^4$ for each fitted scale; parameter values that gave a non-finite error, a solved coefficient
that was not positive, or a solved scale beyond $10^{\pm 6}$ were given an error of $10^3$.  The worst errors
reported are maxima over the check grid, and each frontier model's worst error was confirmed with mpmath at 40 digits.
The whole search takes about eight hours on two processor cores; the notebook gives the settings.

Theorems~\ref{prop:Abounds} and~\ref{prop:Cbounds} are also machine-checked in Lean~4 with Mathlib, building on the
formalization of \citet{lyon2026median}; the Lean project is an ancillary file with the arXiv version of this paper.  The formal
proofs take as hypotheses the published facts the appendices cite (that $\nuG - \alpha$ is nonincreasing and
$\nuG > \alpha - \frac{1}{3}$, that $\nuG$ is convex, and the bounds of \citet{adell2008ramanujan} at integers) and six-digit
brackets on three values, $\nuG(\frac{1}{2}) > 0.227468$, $0.858367 < \nuG'(\frac{1}{2}) < 0.858368$, and $0.968044 < \nuG'(1) < 0.968045$;
everything else, including every numerical check, is proved.  Two
steps take different routes from the appendices.  For Lemma~\ref{lem:akalin}, $\log\Gamma(\alpha + 2)$ is first bounded above by
its Taylor polynomial of degree 4 at $\alpha = 0$, from Euler's product, so that the lemma's argument needs derivatives of
elementary functions only.  For Theorem~\ref{prop:Cbounds} with $\alpha \ge 1$, an upper bound on $2^x - 1 - cx^2$ by a
polynomial in $x = 1/\alpha$, from Taylor's theorem, turns the comparison with the U8 chords into polynomial inequalities, checked
on 36 short intervals, in place of the tangent-line argument.

An exploratory study with Matlab and Google DeepMind's AI co-mathematician \citep{zheng2026ai} compared arctangent and rational interpolators for the
three forms, and reported preliminary Pareto results minimaxing the sum of relative and absolute errors; those results were not used
directly.  In preparing this paper, the author used Claude (Anthropic) to rebuild the computations from scratch as a Python module and
notebook, to run and check the fits, and to draft portions of the text.  Claude found the spurious far-out optima and the rearranged
large-$\alpha$ error expressions.  For Appendix~\ref{app:Abounds}, Claude found the convex-combination identity and the tangent-line
argument, following the author's suggestion to use $\nu_{\Gamma 4}$.  The author directed the work and revised the text, and is
responsible for its content.  A Python notebook that reproduces all fits, tables, figures, and numbers quoted here is included as an
ancillary file with the arXiv version of this paper.  The Python code is entirely as written by Claude in dozens of iterations with
the author, as is the Lean code.

\section{Summary}
\label{sec:summary}

Writing gamma median approximations as interpolations between two elementary functions, each correct in one limit, turns the
problem into approximating an interpolator $g(\alpha)$ that goes from 0 to 1.  For each of the three forms considered, the ideal
interpolator $g^*(\alpha)$ is known in closed form at both ends and at $\alpha = 1$, and interpolators whose coefficients
are all set from what is known, with no fitting and no error criterion,
already reach a worst relative error of $0.0011$ and a worst absolute error of $0.000195$ over all $\alpha > 0$.
Fitting further coefficients to minimize the worst absolute-or-relative error trades accuracy against complexity along a
Pareto frontier, through $0.000087$ with two fitted constants and $0.000066$ with three, to $0.000041$ with five; setting some
coefficients from $g^*$ rather than fitting them gives most of the frontier.  Form C is best with zero and one fitted constants.
With two to five, every frontier model is one of two closely related Form B rationals of degree 4, with no $\alpha^2$ term in
the numerator.  Forms B and C, built on the end function
$\nu_{\Gamma 4}$ from the beta median, beat Form A, the form of \citet{lyon2021closed}, by factors of 6 and 12 with every scale
set, and Form B beats it by factors of 1.8 to 2.5 with two to five fitted constants.  Table~\ref{tab:analytic} lists the analytic approximations, and
Table~\ref{tab:frontier} the frontier approximations with their fitted constants.

For readers who want just a formula, three of the frontier models are as follows; the worst errors are over all $\alpha > 0$.
\begin{center}
\fbox{\begin{minipage}{0.92\linewidth}
\smallskip
\centering
Constants:  $e^{-\gamma} \approx 0.561459$, \ $\log 2 \approx 0.693147$, \
$c = \left(\frac{1}{2}\log 2 - \frac{1}{3}\right)\log 2 \approx 0.0091774$.
\begin{align*}
&\text{no fitted constants, both scales set from $g^*$, Form C, worst error $0.0011$:} \\[2pt]
\nu_\Gamma(\alpha) &\approx \frac{(1 - g)\,e^{-\gamma} + g\log 2}{2^{1/\alpha} - 1 - c/\alpha^2}, \\
g &= \frac{\alpha/0.285173 + (\alpha/0.243693)^2}{1 + \alpha/0.285173 + (\alpha/0.243693)^2} \\[8pt]
&\text{2 fitted constants, 3 scales set from $g^*$, Form B, worst error $0.000087$:} \\[2pt]
\nu_\Gamma(\alpha) &\approx \frac{(1 - g)\,e^{-\gamma} + g\log 2}{2^{1/\alpha} - 1}, \\
g &= \frac{\alpha/0.2851729 + (\alpha/0.1940768)^3 + (\alpha/0.2054328)^4}{1 + (\alpha/0.1549773)^2 + (\alpha/0.2171439)^3 + (\alpha/0.2054328)^4} \\[8pt]
&\text{5 fitted constants, 1 scale set from $g^*$, Form B, worst error $0.000041$:} \\[2pt]
\nu_\Gamma(\alpha) &\approx \frac{(1 - g)\,e^{-\gamma} + g\log 2}{2^{1/\alpha} - 1}, \\
g &= \frac{\alpha/0.285173 + (\alpha/0.180493)^3 + (\alpha/0.196530)^4}{1 + \alpha/6.40642 + (\alpha/0.142651)^2 + (\alpha/0.200966)^3 + (\alpha/0.196530)^4}
\end{align*}

\end{minipage}}
\end{center}
In the first, both scales are set from $g^*$:  $a_1 = 1/g^{*\prime}(0)$ and $a_2 = K^{1/2}$, with
\[
g^{*\prime}(0) = \frac{\pi^2 e^{-\gamma}}{12(\log 2 - e^{-\gamma})} \, , \qquad K = \frac{\log 2}{\log 2 - e^{-\gamma}}\left(\frac{37}{405} - \frac{\log^2 2}{6}\right) \, .
\]
In the second, $a_1$ is set in the same way, and $b_3$ and $b_2$ from Form B's $K$ and $g^*(1) = 1$; in the third, $a_1$
is set in the same way.
For computation, the polynomials are better written in multiplication form, $\alpha/a_1 + (\alpha/a_2)^2 + \cdots =
u_1\alpha + u_2\alpha^2 + \cdots$ with $u_i = a_i^{-i}$, and evaluated by Horner's rule, than by dividing $\alpha$ by each
scale, and $2^{1/\alpha} - 1$ in Forms B and C as $\mathrm{expm1}(\log 2/\alpha)$, which avoids cancellation at
large $\alpha$.  Numerator and denominator then each cost one multiplication and one addition per power of $\alpha$, whether their
scales were fitted, set from $g^*$, or tied, except that the denominator $1 + S$ of $S/(1 + S)$ comes free; so the cost of an
interpolator goes with its degree, and roughly with its total number of scales, not with $n$.  The whole
approximation takes one power $2^{1/\alpha}$, two or three divisions, and a few dozen other operations at most, against the
iterative inversion of the incomplete gamma function by which library routines compute the median.

The approximations also make good starting values for that iterative inversion.  A Newton step on $P(\alpha, x) = \frac{1}{2}$
takes a relative error $\delta$ to about $(1 + \nuG - \alpha)\delta^2/2$, between $\delta^2/3$ and $\delta^2/2$, since
$\alpha - \frac13 < \nuG < \alpha$.  A Halley step costs little more, since $\partial^2 P/\partial x^2 = ((\alpha - 1)/x - 1)\,\partial P/\partial x$,
and reduces $\delta$ to order $\delta^3$.  Over $10^{-3} \le \alpha \le 10^4$, one Halley step from the two-scale analytic
Form C ($0.0011$) leaves a relative error of at most $1.1 \times 10^{-10}$, well within single precision ($6 \times 10^{-8}$),
and two Halley steps reach double precision.  One Newton step from it leaves $5.1 \times 10^{-7}$; from the Form B rational with
two fitted constants on the frontier ($0.000087$), one step leaves $3.2 \times 10^{-9}$, and two reach double precision.
Each step evaluates the incomplete gamma function, at far more cost than any of our formulas, so a more accurate start that
saves a step is worth its extra fitted constants.
Double-precision library routines start instead from more elaborate case-by-case approximations, such as those of
\citet{didonato1986computation}.

\section*{Acknowledgments}

The author thanks Fred Akalin for the proof of Lemma~\ref{lem:akalin}, and acknowledges Claude (Anthropic) and Google DeepMind's AI co-mathematician for help with computation, analysis, proofs, and drafting; their use is described in Section~\ref{sec:numerics}.

\appendix
\section{Form A's ends are bounds}
\label{app:Abounds}

\begin{figure}[b]
\centering
\includegraphics[width=0.9\textwidth]{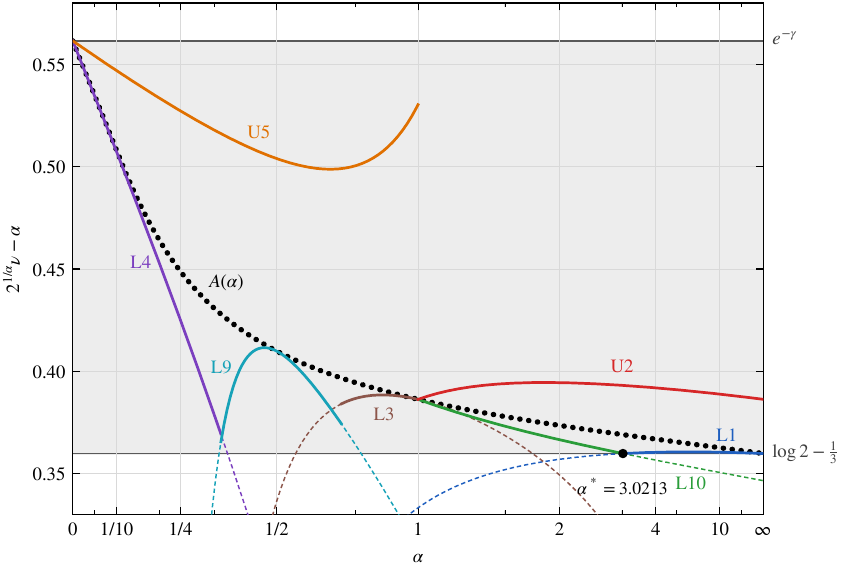}
\caption{The bounds of Table~\ref{tab:abounds} in the box.  Each bound $X$ on $\nuG$ is drawn as $2^{1/\alpha}X - \alpha$, so
that $\nuG$ itself gives $A(\alpha)$ (black dots), and L8 and U6 are the box's floor $\log 2 - \frac{1}{3}$ and ceiling $e^{-\gamma}$.
Lines are thick where the proof uses a bound and dashed where it still holds.  L9 and L3 touch $A$ at $\alpha = \frac{1}{2}$
and 1, and L4 and U5 meet it at $\alpha = 0$.  L10 and L1 cross on the floor, at $\alpha^* = 3.0213$, as the convex
combination in the proof requires.}
\label{fig:abox}
\end{figure}

\citet{lyon2023tight} proved that $A$ lies strictly between its limits, using eight numerically certified point bounds on
$\nuG$ and a published upper bound on the gamma function.  The proof below is shorter.  Besides elementary bounds on the
incomplete gamma function, it uses four known facts:  $\nuG$ is convex \citep{berg2008convexity}; $\nuG' < 1$
\citep{berg2006chen}; $\nuG > \alpha - \frac{1}{3}$ \citep{chen1986bounds}, which also follows from the previous fact and
$\nuG - \alpha \to -\frac{1}{3}$; and Theorem~3 of \citet{lyon2026median}, that $\nuG < \nu_{\Gamma 4}$ for $\alpha < 1$ and
$\nuG > \nu_{\Gamma 4}$ for $\alpha > 1$.  For $\alpha < 1$, the upper bound rests on a lemma proved by the author's colleague
Fred Akalin (Lemma~\ref{lem:akalin}).  The few numbers the proof needs are checked in 50-digit interval arithmetic in the notebook
(\texttt{gamma\_median\_abounds.py}), and the theorem is machine-checked in Lean~4 (Section~\ref{sec:numerics}).  Figure~\ref{fig:abox} shows the bounds in the box between the limits of $A$, and
Table~\ref{tab:abounds} lists them; L9 and L10 are new, and the other names are those of \citet{lyon2023tight}.  (U7 and U8,
used in Appendix~\ref{app:Cbounds}, are also new.)

\begin{table}[tbp]
\caption{Bounds on $\nuG(\alpha)$ used to prove Theorem~\ref{prop:Abounds}.  For $\alpha \ge 1$ the proof uses L1 and
L10 together, whichever is larger; they cross at $\alpha^* = 3.0213$.  The slopes are $\nuG'(\frac{1}{2}) = 0.858367$ and
$\nuG'(1) = 0.968045$, and $\nuG(\frac{1}{2}) = 0.227468$.}
\label{tab:abounds}
\centering
\small
\setlength{\tabcolsep}{4pt}
\renewcommand{\arraystretch}{1.3}
\begin{tabular}{lllll}
\toprule
Name & Bound & Holds for & Used for & From \\
\midrule
L4 & $2^{-1/\alpha}\,\Gamma(1 + \alpha)^{1/\alpha}$ & all $\alpha$ & $\alpha \le 0.35$ & $e^{-t} < 1$ in $P$ \\
L9 & $\nuG(\frac{1}{2}) + \nuG'(\frac{1}{2})(\alpha - \frac{1}{2})$ & all $\alpha$ & $0.35 \le \alpha \le 0.7$ & convexity of $\nuG$ \\
L3 & $\log 2 + \nuG'(1)(\alpha - 1)$ & all $\alpha$ & $0.7 \le \alpha \le 1$ & convexity of $\nuG$ \\
L10 & $\nu_{\Gamma 4} = \log 2/(2^{1/\alpha} - 1)$ & $\alpha > 1$ & $1 \le \alpha \le \alpha^*$ & \citet[Theorem~3]{lyon2026median} \\
L1 & $\alpha - \frac{1}{3}$ & all $\alpha$ & $\alpha \ge \alpha^*$ & \citet{chen1986bounds} \\
L8 & $2^{-1/\alpha}(\log 2 - \frac{1}{3} + \alpha)$ & all $\alpha$ & & Theorem~\ref{prop:Abounds} \\
\midrule
U3 & $\alpha\log 2$ & $\alpha \le 1$ & in U5 & $\nuG < \nu_{\Gamma 4}$, or convexity of $\nuG$ \\
U5 & $2^{-1/\alpha}\,\Gamma(1 + \alpha)^{1/\alpha}\big(1 - \frac{\alpha^2\log 2}{1 + \alpha}\big)^{-1/\alpha}$ & $\alpha \le 1$ & $\alpha \le 1$ & $e^{-t} > 1 - t$ and U3 \\
U2 & $\alpha - 1 + \log 2$ & $\alpha \ge 1$ & $\alpha \ge 1$ & \citet{gaunt2021bounds}; $\nuG' < 1$ \\
U6 & $2^{-1/\alpha}(e^{-\gamma} + \alpha)$ & all $\alpha$ & & Theorem~\ref{prop:Abounds} \\
\bottomrule
\end{tabular}
\end{table}

\begin{theorem}[\citealp{lyon2023tight}]
\label{prop:Abounds}
For all $\alpha > 0$, $\log 2 - \frac{1}{3} < A(\alpha) < e^{-\gamma}$.  Equivalently,
$2^{-1/\alpha}(\log 2 - \frac{1}{3} + \alpha) < \nuG(\alpha) < 2^{-1/\alpha}(e^{-\gamma} + \alpha)$.
\end{theorem}

\begin{proof}
Write $x = 1/\alpha$ and $\ell(\alpha) = 2^{-x}(\Ainf + \alpha)$, so that $A > \Ainf$ exactly where $\nuG > \ell$.

\emph{Lower bound, $\alpha \ge 1$ (L10 and L1).}  Since $(1 - 2^{-x})\,\nu_{\Gamma 4} = 2^{-x}\log 2$,
\[
\ell(\alpha) = (1 - 2^{-x})\,\nu_{\Gamma 4}(\alpha) + 2^{-x}\Big(\alpha - \frac{1}{3}\Big) \, ,
\]
a convex combination of two lower bounds on $\nuG$ for $\alpha > 1$, so $\ell < \nuG$ there.  At $\alpha = 1$,
$\nuG = \log 2 > \frac{1}{2}(\Ainf + 1)$, since $\log 2 > \frac23$.

\emph{Lower bound, $0 < \alpha \le 0.35$ (L4).}  Since $e^{-t} < 1$ in the integral that defines $P$,
$P(\alpha, x) < x^\alpha/\Gamma(\alpha + 1)$, so $\nuG > 2^{-x}\,\Gamma(1 + \alpha)^{1/\alpha}$.  The series
$\log\Gamma(1 + \alpha)/\alpha = -\gamma + \sum_{j \ge 2} (-1)^j\zeta(j)\,\alpha^{j-1}/j$ alternates with decreasing terms for
$0 < \alpha < 1$, so it exceeds its partial sum
$p(\alpha) = -\gamma + \frac{1}{2}\zeta(2)\alpha - \frac{1}{3}\zeta(3)\alpha^2 + \frac14\zeta(4)\alpha^3 - \frac15\zeta(5)\alpha^4$.
On $(0, 0.35]$, $p < -\gamma + \frac{1}{2}\zeta(2)\alpha + \frac14\zeta(4)\alpha^3 < 0$ and
$p' < \frac{1}{2}\zeta(2) + \frac34\zeta(4)\alpha^2 < 1$, so $e^{p} - \Ainf - \alpha$ is strictly decreasing; at $\alpha = 0.35$ it
is $0.00915$.  Hence $\nuG > 2^{-x}e^{p} > \ell$.

\emph{Lower bound, $0.35 \le \alpha \le 1$ (L9 and L3).}  Here we use tangent lines.  Since
\[
\ell'' = 2^{-x}x^4\log 2\,\Big(\Ainf\log 2 - \Big(\log 2 - \frac23\Big)\alpha\Big) \, ,
\]
the function $\ell$ is convex for $\alpha < \Ainf\log 2/(\log 2 - \frac23) = 9.418$.  The median is convex, so it lies above each of its
tangent lines, and a tangent line minus $\ell$ is concave, so it is positive on an interval if it is positive at both ends.
Differentiating $P(\alpha, \nuG(\alpha)) = \frac{1}{2}$ gives
\[
\nuG' = \nuG^{1-\alpha}e^{\nuG}\Big(\frac{1}{2}\,\psi(\alpha)\Gamma(\alpha) - \int_0^{\nuG} t^{\alpha-1}e^{-t}\log t\,dt\Big) \, ,
\]
with $\psi$ the digamma function.  At $\alpha = 1$, where $\nuG = \log 2$ and $\psi(1) = -\gamma$, integrating by parts against
$1 - e^{-t}$ gives
\[
\nuG'(1) = 2\operatorname{Ein}(\log 2) - \log\log 2 - \gamma = 0.968045 \, , \qquad
\operatorname{Ein}(z) = \sum_{n \ge 1} \frac{(-1)^{n+1}z^n}{n \cdot n!} \, ,
\]
or equivalently $\nuG'(1) = \gamma + \log\log 2 + 2E_1(\log 2)$.  At $\alpha = \frac{1}{2}$, $P(\frac{1}{2}, x) = \operatorname{erf}(\sqrt{x})$,
so $\nuG(\frac{1}{2}) = s^2 = 0.227468$ with $s = \operatorname{erf}^{-1}(\frac{1}{2})$.  With $\psi(\frac{1}{2}) = -\gamma - 2\log 2$,
$\Gamma(\frac{1}{2}) = \sqrt\pi$, and $t = u^2$, expanding $e^{-u^2}$ and using $\operatorname{erf}(s) = \frac{1}{2}$ gives
\[
\nuG'\big(\tfrac{1}{2}\big) = s\,e^{s^2}\Big(4sF - \frac{\sqrt\pi}{2}\big(\gamma + \log 4s^2\big)\Big) = 0.858367 \, , \qquad
F = \sum_{n \ge 0} \frac{(-s^2)^n}{n!\,(2n + 1)^2} \, .
\]
The tangent line at $\frac{1}{2}$ exceeds $\ell$ by $0.00075$ at $\alpha = 0.35$ and by $0.0054$ at $0.7$; the tangent line at 1
exceeds it by $0.0090$ at $0.7$ and by $0.0132$ at 1.

\emph{Upper bound, $\alpha \ge 1$ (U2).}  Since $\nuG' < 1$, $\nuG(\alpha) - \alpha \le \nuG(1) - 1 = \log 2 - 1$ (U2, as \citealp{gaunt2021bounds} also
showed), so
$A \le \big(2^x(1 - (1 - \log 2)x) - 1\big)/x$, with $0 < x \le 1$.  Since $e^y \le 1 + y + y^2$ for $0 \le y \le 1$,
$2^x \le 1 + x\log 2 + x^2\log^2 2$; with $\delta = 2\log 2 - 1$ and $y = x\log 2$, this gives
\[
A \le \delta + \delta y - (1 - \log 2)\,y^2 \le \delta + \frac{\delta^2}{4(1 - \log 2)} = 0.507870 < e^{-\gamma} \, .
\]

\emph{Upper bound, $0 < \alpha \le 1$ (U3 and U5).}  At $\alpha = 1$, $A = 2\log 2 - 1 < e^{-\gamma}$.  For $\alpha < 1$,
$\nuG < \nu_{\Gamma 4} = \log 2/(2^x - 1) \le \alpha\log 2$, since $2^x - 1 \ge x$ for $x \ge 1$.  Since $e^{-t} > 1 - t$ in the
integral that defines $P$,
\[
\frac{1}{2} = P(\alpha, \nuG) > \frac{\nuG^\alpha}{\Gamma(\alpha + 1)}\Big(1 - \frac{\alpha\,\nuG}{\alpha + 1}\Big)
> \frac{\nuG^\alpha\,(1 + \alpha - \alpha^2\log 2)}{\Gamma(\alpha + 2)} \, .
\]
By Lemma~\ref{lem:akalin}, $\Gamma(\alpha + 2) < (\alpha + e^{-\gamma})^\alpha(1 + \alpha - \alpha^2\log 2)$, so
$\nuG^\alpha < \frac{1}{2}(\alpha + e^{-\gamma})^\alpha$, that is, $A < e^{-\gamma}$.
\end{proof}

\begin{lemma}[Fred Akalin]
\label{lem:akalin}
With $b = e^{-\gamma}$, the function
$F(\alpha) = \alpha\log(\alpha + b) + \log(1 + \alpha - \alpha^2\log 2) - \log\Gamma(\alpha + 2)$ is positive for
$0 < \alpha \le 1$.
\end{lemma}

\begin{proof}
First, $F(0) = 0$ and $F'(0) = \log b + 1 - \psi(2) = 0$, since $\psi(2) = 1 - \gamma$.  Next,
$F(1) = \log\big((1 + b)(1 - \frac{1}{2}\log 2)\big) > \log(\frac{14}{9}\cdot\frac{13}{20}) = \log\frac{91}{90} > 0$, using
$\frac59 < b < \frac{9}{16}$ and $\log 2 < \frac{7}{10}$.

Factor $1 + \alpha - \alpha^2\log 2 = \log 2\,(\alpha + r)(s - \alpha)$, where $r = 2/(\sqrt{1 + 4\log 2} + 1) = 0.679737 >
\frac{27}{40}$ and $s = 2/(\sqrt{1 + 4\log 2} - 1) = 2.122432$, between 2 and $\frac{15}{7}$.  With
$\psi''(z) = -2\sum_{k \ge 0} (z + k)^{-3}$, the third derivative is $F''' = G - H$, where
\[
G(\alpha) = 2\sum_{k \ge 2} \frac{1}{(\alpha + k)^3} - \frac{2}{(s - \alpha)^3} \, , \qquad
H(\alpha) = \frac{\alpha + 3b}{(\alpha + b)^3} - \frac{2}{(\alpha + r)^3} \, .
\]
The function $G$ is strictly decreasing on $[0, 1]$.  So is $H$:  with $u = \alpha + b$ and $d = r - b > \frac{27}{40} - \frac{9}{16} > \frac19$,
$H' < 0$ is equivalent to $(u + 3b)(1 + d/u)^4 > 3$, and by Bernoulli's inequality and the arithmetic--geometric mean inequality,
\[
(u + 3b)\Big(1 + \frac{d}{u}\Big)^4 \ge (u + 3b)\Big(1 + \frac{4d}{u}\Big) = u + \frac{12bd}{u} + 3b + 4d
\ge 4\sqrt{3bd} + 3b + 4d > \frac49\sqrt{15} + \frac{19}{9} > 3 \, .
\]
So $F''' < G(0) - H(1)$ on $(0, 1)$.  Since $t^{-3}$ is convex, $2\sum_{k \ge 2} k^{-3} < \frac14 + 2\int_{5/2}^\infty t^{-3}\,dt
= \frac{41}{100}$, and $s^3 < 10$, so $G(0) < \frac{41}{100} - \frac15 = \frac{21}{100}$.  Since $(1 + 3t)/(1 + t)^3$ is
decreasing for $t > 0$ and $b < \frac23 < r$, $H(1) > (3 - 2)/(\frac53)^3 = \frac{27}{125} > \frac{21}{100}$.  So $F''' < 0$
on $(0, 1)$.

Finally, $h = \alpha F' - 2F$ has $h(0) = h'(0) = 0$ and $h'' = \alpha F''' < 0$, so $h < 0$ on $(0, 1]$, and
$(F/\alpha^2)' = h/\alpha^3 < 0$.  So $F/\alpha^2$ is strictly decreasing on $(0, 1]$, and $F(\alpha) > \alpha^2 F(1) > 0$.
\end{proof}

\paragraph{Corrections to \citet{lyon2023tight}.}  We take this opportunity to correct some errors in that paper, none of
which affects its theorems or their proofs; it writes $k$ for the shape parameter.  The last sentence before its section on
recent related work should name Theorems U6 and L8, not U6 and U8, and Theorem U6 should say that the median is bounded
above by $U_6(k)$, not $L_6(k)$.  In the statement of Theorem L3, the slope is
$\nu'(1) = \gamma - 2\operatorname{Ei}(-\log 2) + \log\log 2 \approx 0.968045$, as at the end of its proof, not
$\gamma - 2\operatorname{Ei}(-\log 2) - \log\log 2$.  In the proof of Theorem L6, the subtracted term after the choice of
$L^b(k)$ has $L^a(k)^k$, not $U^a(k)^k$; and in the formulas for $L_6(k)$ at the end of that proof and for $L_{61}(k)$ in
Corollary~1, the factor $k^2/(k + 1)$ should be $(1 - e^{-k})/(k + 1)$, as in the statement of Theorem L6 and in Table~1
there.  In the proof of Lemma~3, the floor of the box is $A_L = \log 2 - \frac{1}{3}$, not $L_U$; in the proof of Lemma~4,
$A_{L4}(k)$ approaches $e^{-\gamma}$, not $e^{\gamma}$, as $k \to 0$; and the caption of Figure~3 there should refer to
Theorem L7, not Theorem 7.

\section{Form C's ends are bounds}
\label{app:Cbounds}

\begin{figure}[b]
\centering
\includegraphics[width=0.9\textwidth]{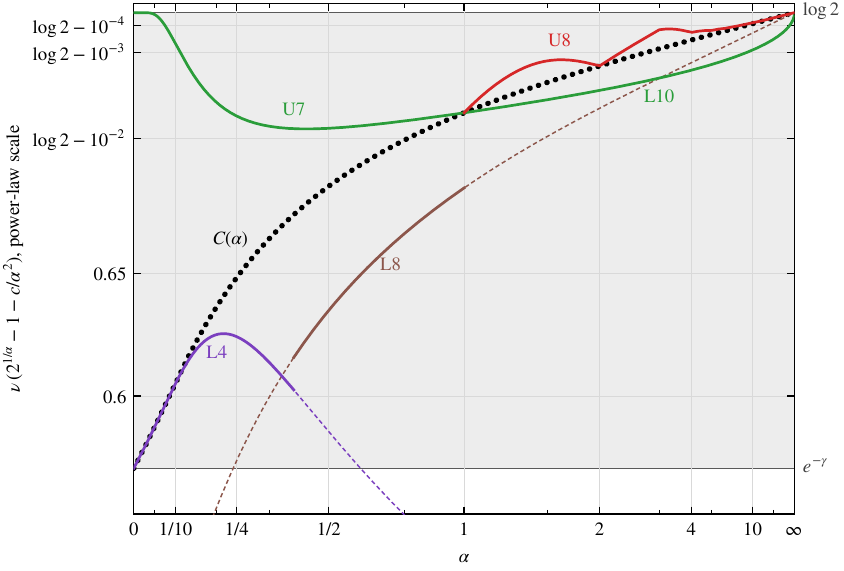}
\caption{The bounds of Table~\ref{tab:cbounds} in the box.  Each bound $X$ on $\nuG$ is drawn as $X(2^{1/\alpha} - 1 - c/\alpha^2)$,
so that $\nuG$ itself gives $C(\alpha)$ (black dots), on a power-law scale, $1 - \big((\log 2 - C)/(\log 2 - e^{-\gamma})\big)^{1/2}$,
that puts the floor $e^{-\gamma}$ at 0 and the ceiling $\log 2$ at 1.  Since $\log 2 - C = O(\alpha^{-2})$, and
$\frac{\pi}{2} - \arctan\alpha \sim 1/\alpha$, $C$ approaches the upper right corner along a straight line.  Lines are thick
where the proof uses a bound and dashed where it still holds.  L10 and U7, both $\nu_{\Gamma 4}$, meet $C$ at $\alpha = 1$; the
chords U8 meet the upper bounds on $\nuG$ at the integers.}
\label{fig:cbox}
\end{figure}

Figure~\ref{fig:cbox} shows the bounds this proof uses, mapped into the box between the limits of $C$, and
Table~\ref{tab:cbounds} lists them.  L4, L8, and L10 are as in Table~\ref{tab:abounds}; U7 is the same function as L10,
$\nu_{\Gamma 4}$, on the other side of $\alpha = 1$, and U8 is new.  Like Theorem~\ref{prop:Abounds}, the theorem is
machine-checked in Lean~4 (Section~\ref{sec:numerics}).

\begin{table}[hbp]
\caption{Bounds on $\nuG(\alpha)$ used to prove Theorem~\ref{prop:Cbounds}.  U8 is the chord of $\nuG$ on $[k, k + 1]$,
$k = 1, 2, \ldots$, through upper bounds $\bar\nu_k$ on $\nuG(k)$:  $\bar\nu_1 = \log 2$, $\bar\nu_2 = 1.6783469900167$, and
$\bar\nu_k = k - \frac{1}{3} + \frac{8}{405(k - 1)}$ for $k \ge 3$.}
\label{tab:cbounds}
\centering
\small
\renewcommand{\arraystretch}{1.3}
\begin{tabular}{llll}
\toprule
Name & Bound & Used for & From \\
\midrule
L4 & $2^{-1/\alpha}\,\Gamma(1 + \alpha)^{1/\alpha}$ & $\alpha \le 0.4$ & $e^{-t} < 1$ in $P$ \\
L8 & $2^{-1/\alpha}(\log 2 - \frac{1}{3} + \alpha)$ & $0.4 \le \alpha \le 1$ & Theorem~\ref{prop:Abounds} \\
L10 & $\nu_{\Gamma 4} = \log 2/(2^{1/\alpha} - 1)$ & $\alpha \ge 1$ & \citet[Theorem~3]{lyon2026median} \\
\midrule
U7 & $\nu_{\Gamma 4} = \log 2/(2^{1/\alpha} - 1)$ & $\alpha \le 1$ & \citet[Theorem~3]{lyon2026median} \\
U8 & $\bar\nu_k + (\bar\nu_{k+1} - \bar\nu_k)(\alpha - k)$ on $[k, k + 1]$ & $\alpha \ge 1$ & convexity of $\nuG$;
\citet{adell2008ramanujan} \\
\bottomrule
\end{tabular}
\end{table}

\begin{theorem}
\label{prop:Cbounds}
For all $\alpha > 0$, $e^{-\gamma} < C(\alpha) < \log 2$.  Equivalently, Form C's end functions
$e^{-\gamma}/(2^{1/\alpha} - 1 - c/\alpha^2)$ and $\log 2/(2^{1/\alpha} - 1 - c/\alpha^2)$ are strict lower and upper bounds on
$\nuG(\alpha)$.
\end{theorem}

\begin{proof}
Write $x = 1/\alpha$ and $D(\alpha) = 1 - 2^{-x}(1 + cx^2)$, so that $C = 2^x\nuG\,D$; $D > 0$ since the denominator in
\eqref{eq:ABC} is positive.  Each step below uses a known bound on $\nuG$; the few numbers it needs are elementary functions at
rational points, which the notebook evaluates in 50-digit interval arithmetic (\texttt{gamma\_median\_cbounds.py}).

\emph{Lower bound, $0 < \alpha \le 0.4$ (L4).}  As in Appendix~\ref{app:Abounds}, $\nuG > 2^{-x}\,\Gamma(1 + \alpha)^{1/\alpha}$, so
$C > \Gamma(1 + \alpha)^{1/\alpha} D$.  The series $\log\Gamma(1 + \alpha)/\alpha + \gamma = \sum_{j \ge 2} (-1)^j\zeta(j)\,\alpha^{j-1}/j$
alternates with decreasing terms for $0 < \alpha < 1$, so $\Gamma(1 + \alpha)^{1/\alpha} \ge
e^{-\gamma}\exp(\frac{\pi^2}{12}\alpha - \frac{1}{3}\zeta(3)\alpha^2) \ge e^{-\gamma}(1 + \kappa\alpha)$ with
$\kappa = \frac{\pi^2}{12} - 0.4\,\zeta(3)/3 = 0.662193$.  Also $1 - D = \alpha f(x)$ with $f(t) = 2^{-t}(t + ct^3)$, and
$2^t f'(t) = -(t\log 2 + ct^3\log 2 - 1 - 3ct^2) < 0$ for $t \ge 2.5$, since the bracket is $0.660$ at $t = 2.5$ and has
derivative $\log 2 + 3ct(t\log 2 - 2) > 0.67$ there.  So $1 - D \le f(2.5)\,\alpha = 0.467291\,\alpha <
\kappa\alpha/(1 + \kappa\alpha)$, the last because $\kappa/(1 + 0.4\kappa) = 0.523523$; hence $(1 + \kappa\alpha)D > 1$ and
$C > e^{-\gamma}$.

\emph{Lower bound, $0.4 \le \alpha \le 1$ (L8).}  By Theorem~\ref{prop:Abounds}, $A > \Ainf$, so $C = (A + \alpha)D > (\Ainf + \alpha)D$.
$D$ is strictly decreasing, since $\frac{d}{d\alpha}\,2^{-x}(1 + cx^2) = 2^{-x}x^2(\log 2 - 2cx + cx^2\log 2) > 0$, the quadratic
in $x$ having discriminant $4c(c - \log^2 2) < 0$.  So on $[s, t]$, $C > (\Ainf + s)\,D(t)$, which exceeds $e^{-\gamma}$ on each
of the seven intervals between $0.4$, $0.45$, $0.5$, $0.6$, $0.7$, $0.8$, $0.9$, and $1$, by at least $0.0206$.

\emph{Lower bound, $\alpha \ge 1$ (L10).}  Here $\nuG \ge \nu_{\Gamma 4}$ \citep[Theorem~3]{lyon2026median}, so
$C \ge \log 2\,\big(1 - cx^2/(2^x - 1)\big) \ge \log 2\,\big(1 - c/(\log 2 + \frac{1}{2}\log^2 2)\big) = 0.686332 > e^{-\gamma}$,
using $(2^x - 1)/x^2 \ge \log 2/x + \frac{1}{2}\log^2 2$ and $x \le 1$.

\emph{Upper bound, $\alpha \le 1$ (U7).}  Here $\nuG \le \nu_{\Gamma 4}$ \citep[Theorem~3]{lyon2026median}, so
$C < B \le \log 2$.

\emph{Upper bound, $1 \le \alpha \le 10$ (U8).}  The upper end function $U(\alpha) = \log 2/(2^x - 1 - cx^2)$ is convex for
$\alpha \ge 1$:  $U = \alpha/q(x)$ with $q(x) = 1 + x/3 + \sum_{j \ge 3} (x\log 2)^{j-1}/j!$, and $U'' = x^3\psi''(x)$ with
$\psi = 1/q$ and $\psi'' = (2q'^2 - qq'')/q^3$; since $q$, $q'$, and $q''$ increase, $2q'(a)^2 - q(b)q''(b) > 0$ on ten intervals
$[a, b]$ of width $0.1$ covering $[0, 1]$ shows that $\psi'' > 0$.  The median is strictly convex \citep{berg2008convexity}, so on
$[k, k + 1]$ it lies below the chord through upper bounds on its values at the ends:  $\nuG(1) = \log 2$;
$\nuG(2) < 1.6783469900167$, from $(1 + x)e^{-x} = \frac{1}{2}$ (Section~\ref{sec:limits}); and
$\nuG(k) < k - \frac{1}{3} + \frac{8}{405(k - 1)}$ for integers $k \ge 3$ \citep{adell2008ramanujan,adell2009inequalities}.  On each
$[k, k + 1]$, $U$ minus the chord is convex, hence at least its tangent line at a point near its minimum; for $k = 1$, \dots, 9 that
tangent line stays above $0.0006$ on the interval.

\emph{Upper bound, $\alpha \ge 10$ (U8).}  With $q = 1 + u$, $u = x/3 + R(x)$, $R \ge 0$, and $1/(1 + u) \ge 1 - u + u^2 - u^3$,
\[
U - \Big(\alpha - \frac{1}{3}\Big) \ge x\Big(\frac19 - \frac{R}{x^2} - x\Big(\frac{1}{3} + \frac{R}{x}\Big)^3\Big) \ge 0.025566\,x
\quad \text{for } x \le 0.1 \, ,
\]
using $R/x^2 \le \frac16\log^2 2 + \frac{1}{24}x\log^3 2\;e^{x\log 2}$.  On $[k, k + 1]$ with $k \ge 10$, the chord of the
previous step is at most $\alpha - \frac{1}{3} + \frac{8}{405(k - 1)}$, and $0.025566/(k + 1) > \frac{8}{405(k - 1)}$.
\end{proof}

\begin{remark}
\label{rem:AC}
Since $C = (A + \alpha)D$ with $D > 0$ strictly decreasing, $C' = (A' + 1)D + (A + \alpha)D'$.  So $C' > 0$ exactly where
$A' > -1 + (A + \alpha)|D'|/D$, a lower bound on $A'$, while the conjecture that $A$ is strictly decreasing is the
upper bound $A' < 0$.  The two statements bound $A'$ from opposite sides, and neither follows from the other this way; strict
increase of $C$ would, however, imply that $A + \alpha = 2^{1/\alpha}\nuG$ is strictly increasing.
\end{remark}

\bibliographystyle{plainnat}
\bibliography{gamma_median}

\begin{thebibliography}{22}
\providecommand{\natexlab}[1]{#1}
\providecommand{\url}[1]{\texttt{#1}}
\expandafter\ifx\csname urlstyle\endcsname\relax
  \providecommand{\doi}[1]{doi: #1}\else
  \providecommand{\doi}{doi: \begingroup \urlstyle{rm}\Url}\fi

\bibitem[Adell and Jodr{\'a}(2008)]{adell2008ramanujan}
J~Adell and P~Jodr{\'a}.
\newblock On a {Ramanujan} equation connected with the median of the gamma
  distribution.
\newblock \emph{Transactions of the American Mathematical Society},
  360\penalty0 (7):\penalty0 3631--3644, 2008.

\bibitem[Adell and Alzer(2009)]{adell2009inequalities}
Jos{\'e}~A. Adell and Horst Alzer.
\newblock Inequalities for the median of the gamma distribution.
\newblock \emph{Journal of Computational and Applied Mathematics}, 232\penalty0
  (2):\penalty0 481--495, 2009.
\newblock \doi{10.1016/j.cam.2009.06.024}.

\bibitem[Berg and Pedersen(2006)]{berg2006chen}
Christian Berg and Henrik~L. Pedersen.
\newblock The {Chen--Rubin} conjecture in a continuous setting.
\newblock \emph{Methods and Applications of Analysis}, 13\penalty0
  (1):\penalty0 63--88, 2006.
\newblock \doi{10.4310/MAA.2006.v13.n1.a4}.

\bibitem[Berg and Pedersen(2008)]{berg2008convexity}
Christian Berg and Henrik~L. Pedersen.
\newblock Convexity of the median in the gamma distribution.
\newblock \emph{Arkiv f{\"o}r Matematik}, 46\penalty0 (1):\penalty0 1--6, 2008.

\bibitem[Chen(2017)]{chen2017median}
Chao-Ping Chen.
\newblock The median of gamma distribution and a related {Ramanujan} sequence.
\newblock \emph{The Ramanujan Journal}, 44\penalty0 (1):\penalty0 75--88, 2017.

\bibitem[Chen and Rubin(1986)]{chen1986bounds}
Jeesen Chen and Herman Rubin.
\newblock Bounds for the difference between median and mean of gamma and
  {Poisson} distributions.
\newblock \emph{Statistics \& Probability Letters}, 4\penalty0 (6):\penalty0
  281--283, 1986.
\newblock \doi{10.1016/0167-7152(86)90044-1}.

\bibitem[Choi(1994)]{choi1994medians}
Kwok~Pui Choi.
\newblock On the medians of gamma distributions and an equation of {Ramanujan}.
\newblock \emph{Proceedings of the American Mathematical Society}, 121\penalty0
  (1):\penalty0 245--251, 1994.
\newblock \doi{10.1090/S0002-9939-1994-1195477-8}.

\bibitem[Corless et~al.(1996)Corless, Gonnet, Hare, Jeffrey, and
  Knuth]{corless1996lambert}
R.~M. Corless, G.~H. Gonnet, D.~E.~G. Hare, D.~J. Jeffrey, and D.~E. Knuth.
\newblock On the {Lambert} {W} function.
\newblock \emph{Advances in Computational Mathematics}, 5\penalty0
  (1):\penalty0 329--359, 1996.
\newblock \doi{10.1007/BF02124750}.

\bibitem[DiDonato and Morris(1986)]{didonato1986computation}
Armido~R DiDonato and Alfred~H Morris, Jr.
\newblock Computation of the incomplete gamma function ratios and their
  inverse.
\newblock \emph{ACM Transactions on Mathematical Software}, 12\penalty0
  (4):\penalty0 377--393, 1986.
\newblock \doi{10.1145/22721.23109}.

\bibitem[DLMF()]{dlmf}
DLMF.
\newblock {NIST} digital library of mathematical functions.
\newblock \url{https://dlmf.nist.gov/}, 2026.
\newblock F.~W.~J. Olver et al., eds.; Release 1.2.8 of 2026-09-15.

\bibitem[Gaunt and Merkle(2021)]{gaunt2021bounds}
Robert~E Gaunt and Milan Merkle.
\newblock On bounds for the mode and median of the generalized hyperbolic and
  related distributions.
\newblock \emph{Journal of Mathematical Analysis and Applications},
  493\penalty0 (1):\penalty0 124508, 2021.

\bibitem[Harris et~al.(2020)Harris, Millman, van~der Walt,
  et~al.]{harris2020array}
Charles~R. Harris, K.~Jarrod Millman, St{\'e}fan~J. van~der Walt, et~al.
\newblock Array programming with {NumPy}.
\newblock \emph{Nature}, 585:\penalty0 357--362, 2020.
\newblock \doi{10.1038/s41586-020-2649-2}.

\bibitem[Lyon(2021)]{lyon2021closed}
Richard~F Lyon.
\newblock On closed-form tight bounds and approximations for the median of a
  gamma distribution.
\newblock \emph{PLOS One}, 16\penalty0 (5):\penalty0 e0251626, 2021.
\newblock \doi{10.1371/journal.pone.0251626}.

\bibitem[Lyon(2023)]{lyon2023tight}
Richard~F Lyon.
\newblock Tight bounds for the median of a gamma distribution.
\newblock \emph{{PLOS One}}, 18\penalty0 (9):\penalty0 e0288601, 2023.
\newblock \doi{10.1371/journal.pone.0288601}.

\bibitem[Lyon(2026)]{lyon2026median}
Richard~F. Lyon.
\newblock On the median of the beta distribution in relation to the median of
  the gamma distribution.
\newblock arXiv:2610.39202, 2026.

\bibitem[Nelder and Mead(1965)]{nelder1965simplex}
John~A. Nelder and Roger Mead.
\newblock A simplex method for function minimization.
\newblock \emph{The Computer Journal}, 7\penalty0 (4):\penalty0 308--313, 1965.
\newblock \doi{10.1093/comjnl/7.4.308}.

\bibitem[Storn and Price(1997)]{storn1997differential}
Rainer Storn and Kenneth Price.
\newblock Differential evolution -- a simple and efficient heuristic for global
  optimization over continuous spaces.
\newblock \emph{Journal of Global Optimization}, 11\penalty0 (4):\penalty0
  341--359, 1997.
\newblock \doi{10.1023/A:1008202821328}.

\bibitem[{The mpmath development team}(2026)]{mpmath}
{The mpmath development team}.
\newblock mpmath: a {P}ython library for arbitrary-precision floating-point
  arithmetic, 2026.
\newblock URL \url{https://mpmath.org/}.
\newblock Version 1.4.0.

\bibitem[Trefethen(2019)]{trefethen2019approximation}
Lloyd~N. Trefethen.
\newblock \emph{Approximation Theory and Approximation Practice, Extended
  Edition}.
\newblock SIAM, Philadelphia, 2019.
\newblock \doi{10.1137/1.9781611975949}.

\bibitem[Virtanen et~al.(2020)Virtanen, Gommers, Oliphant,
  et~al.]{virtanen2020scipy}
Pauli Virtanen, Ralf Gommers, Travis~E. Oliphant, et~al.
\newblock {SciPy} 1.0: fundamental algorithms for scientific computing in
  {Python}.
\newblock \emph{Nature Methods}, 17:\penalty0 261--272, 2020.
\newblock \doi{10.1038/s41592-019-0686-2}.

\bibitem[You(2017)]{you2017approximation}
Xu~You.
\newblock Approximation of the median of the gamma distribution.
\newblock \emph{Journal of Number Theory}, 174:\penalty0 487--493, 2017.

\bibitem[Zheng et~al.(2026)Zheng, von Glehn, Zwols, Beloshapka, Buesing, Roy,
  Wattenberg, Georgiev, Schmidt, Cowie, et~al.]{zheng2026ai}
Daniel Zheng, Ingrid von Glehn, Yori Zwols, Iuliya Beloshapka, Lars Buesing,
  Daniel~M. Roy, Martin Wattenberg, Bogdan Georgiev, Tatiana Schmidt, Andrew
  Cowie, et~al.
\newblock {AI} co-mathematician: Accelerating mathematicians with agentic {AI}.
\newblock arXiv:2605.06651, 2026.

\end{thebibliography}

\end{document}